\documentclass[a4paper,11pt,twoside]{article}
\usepackage{amssymb,amsmath,amsthm}
\usepackage{caption}
\usepackage{subcaption}
\usepackage{arydshln}
\usepackage{tikz}
\usetikzlibrary{shapes}
\usepackage[a4paper,includeheadfoot,margin=2cm]{geometry}

\newtheorem{proposition}{\bf Proposition}
\newtheorem{remark}{\bf Remark}
\def\R{\mathbb R}
\newcommand{\pr}[1]{{\rm r}^{(#1)}} 
\newcommand{\pj}[1]{{\rm j}^{(#1)}} 
\def\pri{{\rm r}} 
\def\pji{{\rm j}} 
\title{An Asymptotic Preserving Two-Dimensional Staggered Grid Method for multiscale transport equations} 

\author{Kerstin K\"upper\thanks{
Department of Mathematics
and Center for Computational Engineering Science,
Schinkelstrasse 2,
D-52062 Aachen, Germany.
(kuepper@mathcces.rwth-aachen.de, frank@mathcces.rwth-aachen.de). The authors were funded by the Excellence Initiative of the German Federal and State Governments.}
 \and Martin Frank\footnotemark[1]
 \and Shi Jin\thanks{ 
        Department of Mathematics,  University of Wisconsin,
                Madision, WI 53706, USA and Department of Mathematics, 
Institute of Natural 
Sciences and MOE-LSC,
Shanghai Jiao Tong University, Shanghai 200240, China (jin@math.wisc.edu). The research of this author was partially supported by NSF DMS grant 1107291: RNMS ``KI-Net'', and NSFC grant 91330203.}}

\date{}
\begin{document}
\maketitle

\begin{abstract}
We propose a two-dimensional asymptotic preserving scheme for linear transport
equations with diffusive scalings. It is an extension of the time
splitting developed  by Jin, Pareschi and Toscani~\cite{Jin-Pareschi-Toscani-2000}, 
but uses spatial discretizations on staggered grids, which preserves
the discrete diffusion limit with a more compact stencil.
The first novelty of this paper is that we propose a staggering in two dimensions that requires fewer unknowns than one could have
naively expected. The second contribution of this paper is that we rigorously analyze the scheme of Jin, Pareschi, and Toscani \cite{Jin-Pareschi-Toscani-2000}. We show that the scheme is AP and obtain an explicit CFL condition, which couples a hyperbolic and a parabolic condition. This type of condition is common for asymptotic preserving schemes and guarantees uniform stability with respect to the mean free path. 
In addition, we obtain an upper bound on the relaxation parameter, which is the crucial parameter of the used time discretization. Several numerical examples
are provided to verify the accuracy and asymptotic property of the scheme.
\end{abstract}

\pagestyle{myheadings}
\thispagestyle{plain}
\markboth{KERSTIN K\"UPPER, MARTIN FRANK, AND SHI JIN}{ASYMPTOTIC PRESERVING TWO-DIMENSIONAL STAGGERED GRID METHOD}

\section{Introduction}
The linear transport equation models particles interacting with a background medium (e.g.,\ neutron transport, linear radiative transfer). In general, the model in scaled variables can be written as~\cite{Jin-Pareschi-Toscani-2000}
\begin{equation} \label{eq:transport}
 \varepsilon\partial_t f + v\cdot\nabla_\mathbf{x} f = \frac 1 \varepsilon \left[ \frac {\sigma_s}{2\pi} \int_\Omega f dv' - \sigma_t f \right] + \varepsilon Q,
\end{equation}
where $f(t,\mathbf{x},v)$ denotes the probability density distribution depending on time $t$, position $\mathbf{x}\in\R^2$, and direction of velocity $v=(\xi,\eta)\in\Omega = \{(\xi,\eta):\; -1\leq \xi,\eta\leq1,\;\xi^2+\eta^2=1\}$ (two-dimensional flatland model -- the extension to three dimensions is straightforward). Moreover, $\sigma_t=\sigma_s+\varepsilon^2 \sigma_a$ is the total transport coefficient, $\sigma_s$ is the scattering coefficient, $\sigma_a$ is the absorption coefficient, and $Q$ is a $v$-independent source term. It is well-known that the limiting equation ($\varepsilon\to 0$) of~\eqref{eq:transport} is the diffusion equation:
\begin{equation}\label{eq:diffusion}
\partial_t \rho = \tfrac12 \nabla_\mathbf{x} \cdot \left(\frac{1}{\sigma_t} \nabla_\mathbf{x} \rho\right) - \sigma_a \rho + Q,
\end{equation}
where $\rho(t,\mathbf{x}) = \frac1{2\pi}\int_\Omega f(t,\mathbf{x},v) dv$. 

In many applications, the scaling parameter of the transport equation $\varepsilon$ (mean free path) may differ in several orders of magnitude, ranging from the rarefied kinetic regime to the hydrodynamic diffusive regime. When $\varepsilon$ is small, in the diffusive regime, the equation becomes numerically
stiff, which leads to numerical challenges: straightforward explicit implementations lead to high computational costs in the diffusive regimes; fully implicit schemes could be difficult to implement 
\cite{Deng-2014}; multiscale multiphysics domain-decomposition approaches, which couple models at different scales, have difficulties in the transition zones, since they need to
transfer data from one scale to another \cite{Jin-2012}. Thus, it is desirable to develop schemes which are suitable for all regimes (no domain-decomposition), but do not require a resolved grid in space and small time compared to the mean free path. This is the objective of asymptotic preserving (AP) schemes.

A scheme is called AP if it preserves the discrete analogue of the asymptotic
transition from the microscopic scale to the mascroscopic one \cite{Jin-1999, Jin-2012}, 
namely,  in the limit $\varepsilon \to 0$, the  discretization of the above transport equation~\eqref{eq:transport} should yield a discretization of the diffusion equation~\eqref{eq:diffusion}.  Such schemes allow mesh sizes and time
steps much bigger than the mean free path/time, yet still capture the correct
physical behavior. The development of such schemes started with steady problems
of linear transport equations by Larsen, Morel, and Miller~\cite{Larsen-Morel-Miller-1987} and for boundary value problems by
Jin and Levermore \cite{Jin-Levermore-1991, Jin-Levermore-1993}. Uniform convergence with respect to the 
mean free path for an AP scheme was first established by Golse, Jin, and Levermore
\cite{Golse-Jin-Levermore-1999}. 

In \cite{Larsen-Morel-Miller-1987}, and its follow-up \cite{Larsen-Morel-1989}, several space discretizations for steady transport problems were investigated, among them diamond differencing and a linear discontinuous Galerkin (LD) method, both of which were identified to behave well in the asymptotic regime. Furthermore, LD leads to a compact stencil for the diffusion equation. 

For time-dependent problems,  AP schemes 
were first designed for nonlinear hyperbolic systems with relaxation by
Jin and Levermore \cite{Jin-1995,Jin-Levermore-1996}. There one needs to
design both the time and the spatial discretization carefully \cite{Jin-Levermore-1996}, in
particular, to overcome the stiffness of the source term.  AP schemes for
time-dependent transport equations with diffusive scaling started by
Jin, Pareschi, and Toscani \cite{Jin-Pareschi-Toscani-1998, Jin-Pareschi-Toscani-2000} and Klar \cite{Klar-1998}.
Since then there have been many new developments in the construction of
AP schemes for a large class of kinetic equations (cf.\ reviews by Degond~\cite{Degond-2013} and Jin~\cite{Jin-2012}). Time discretizations usually need 
an implicit-explicit (IMEX) approach \cite{Caflisch-Jin-Russo-1997, Jin-1995, Pareschi-Russo-2000,Boscarino-Russo-2009}, exponential integration methods~\cite{Dimarco-Pareschi-2011}, BGK-type penalty methods~\cite{Filbet-Jin-2010}, or micro-macro decomposition-based schemes~\cite{Lemou-Mieussens-2008,Liu-Mieussens-2010}. See also \cite{Mieussens-2013, Xu-Huang-2010}.  One key idea of the schemes is to split the
equation into a nonstiff part, which is treated explicitly, and a stiff part,
which will be implicit but can be implemented explicitly. The splitting
should be taken in a way such that the combination preserves the AP property. Another possibility to treat the time-dependent case, which has been used often in the neutron transport literature, is to use an implicit semidiscretization in time first (e.g.,\ backward Euler). This reduces the problem to the successive solution of steady problems (with a modified absorption term). For these, one can use a variety of techniques, among them the second order self-adjoint form of the transport equation. For a recent example, see \cite{Morel-Adams-Noh-McGhee-Evans-Urbatsch-2006}, where different spatial discretizations are investigated. The second order self-adjoint form, however, does not exist for the fully time-dependent, undiscretized transport equation.

In this paper, we present a two-dimensional AP scheme for the time-dependent transport equation~\eqref{eq:transport}, where we do not use a semidiscretization in time. We combine a well-known scheme~\cite{Jin-Pareschi-Toscani-2000} which is based on the parity equations (which in turn are well-known~\cite{Lewis-Miller-1984}) with staggered grids, which in one dimension are fully understood~\cite{Jin-2012}. As we have mentioned above, the advantage of the staggered approach, compared  with the regular grid approach in \cite{Jin-Pareschi-Toscani-2000}, is that in the diffusion limit we approach a {\it compact} stencil, as pointed out by Jin \cite{Jin-2012}.
To be more precise, in one space dimension, using a regular Cartesian grid, the discrete diffusion limit of the scheme of \cite{Jin-Pareschi-Toscani-2000}
approximates the diffusion operator in  (\ref{eq:diffusion}) by $(\rho_{i+2}-2\rho_i+\rho_{i-2})/(2\Delta x)^2$ (for the case $\sigma_t\equiv 1 $), while the
current scheme gives the compact discretization  $(\rho_{i+1}-2\rho_i+\rho_{i-1})/(\Delta x)^2$ which offers a better spatial resolution.

However, it is not trivial to extend the scheme from \cite{Jin-Pareschi-Toscani-2000} to staggered grids in two dimensions. The first novelty of this paper is that we propose a staggering in two dimensions that requires fewer unknowns than one could have
naively expected. The scheme shares similarities with diamond differencing but turns out to be different. The second contribution of this paper is that we present a rigorous stability analysis of the scheme \cite{Jin-Pareschi-Toscani-2000} in one dimension and show uniform stability. Previously, only a stability argument was available. Similar to~\cite{Lemou-Mieussens-2008,Liu-Mieussens-2010}, we obtain an explicit CFL condition, which couples a hyperbolic and a parabolic condition and guarantees uniform stability. The scheme contains a splitting parameter that distributes terms between the explicit and implicit parts. For the choice of this parameter, our stability analysis yields an upper bound that is more restrictive than the one that has previously been used~\cite{Jin-Pareschi-Toscani-2000}. We discuss the choices of the methods we combine in the text. The advantages and disadvantages of these have been described in the original papers. For the sake of completeness, we repeat them in the text.

The remainder of the paper is organized as follows. In section~\ref{sec:method}, we first derive the parity equations for the linear transport equation. Then, we describe the numerical method. In section~\ref{sec:AP-property}, we show the AP property. We consider the asymptotic limit of the scheme and we state and prove a stability result, which gives a CFL condition and an upper bound on the relaxation parameter. Finally, several numerical tests are presented in section~\ref{sec:numerics} to confirm the AP property of the scheme. Before the conclusion, we comment on several possible extensions (section \ref{sec:ext}).

\section{The Numerical Method}
\label{sec:method}
First, we reformulate the transport equation into a parity equation. Then, we describe the angular, the spatial, and the time discretization of the method in detail.

We begin with the transport equation in the diffusive scaling~\eqref{eq:transport}, restricted to two spatial dimensions $\mathbf{x} = (x,y)$. Following \cite{Jin-Pareschi-Toscani-2000}, the equation is split into four parts according to the quadrants of the velocity space. We obtain four equations with non-negative $\xi,\eta\geq0$. This system can be rewritten if we define the even and odd parities
\begin{equation}\label{eq:definition-parities}
\begin{array}{ll}
\pr{1}(\xi, \eta)= \frac 1 2 [ f(\xi,-\eta)+ f(-\xi,\eta)],& \pr{2}(\xi, \eta)= \frac 1 2 [ f(\xi,\eta)+ f(-\xi,-\eta)],\\[1ex]
\pj{1}(\xi, \eta)= \frac 1 {2\varepsilon} [ f(\xi,-\eta)- f(-\xi,\eta)],& \pj{2}(\xi, \eta)= \frac 1 {2\varepsilon} [ f(\xi,\eta)-f(-\xi,-\eta)],
\end{array}
\end{equation}
leading to
\begin{subequations} \label{eq:pde_system}
\begin{align}
\partial_t \pr{1} + \xi \partial_x \pj{1} -\eta \partial_y \pj{1} &= -\frac{\sigma_s}{\varepsilon^2} ( \pr{1} - \rho)-\sigma_a \pr{1} +Q,\label{eq:r1}\\
\partial_t \pr{2} + \xi \partial_x \pj{2} +\eta \partial_y \pj{2} &= -\frac{\sigma_s}{\varepsilon^2} (\pr{2} - \rho)-\sigma_a \pr{2} +Q,\label{eq:r2}\\
\partial_t \pj{1} + \frac\xi{\varepsilon^2} \partial_x \pr{1} -\frac\eta{\varepsilon^2} \partial_y \pr{1} &= -\frac{\sigma_s}{\varepsilon^2} \pj{1}-\sigma_a \pj{1},\label{eq:j1}\\
\partial_t \pj{2} + \frac\xi{\varepsilon^2} \partial_x \pr{2} +\frac\eta{\varepsilon^2} \partial_y \pr{2} &= -\frac{\sigma_s}{\varepsilon^2} \pj{2} -\sigma_a  \pj{2},\label{eq:j2}
\end{align}
\end{subequations}
where $\rho= \tfrac 1{2\pi} \int_{|v|=1} f dv.$ It may seem that we have quadrupled the number of unknowns.
Note, however, that due to symmetry we need to solve these equations for $\xi,\eta$ in the positive quadrant only. Thus the number of unknowns is effectively the same.

The discretization of the angular variable uses Gaussian quadrature points. Let
\begin{equation}
\xi(\lambda)=\cos(\lambda \pi/2) \quad\text{and}\quad \eta(\lambda)=\sin(\lambda\pi/2) \label{eq:lambda}
\end{equation}
for every $0\leq\lambda\leq1$. Then, the density
\begin{equation}
\rho= \tfrac12 \int_0^1 [\pr{1}(\xi,\eta)+\pr{2}(\xi,\eta)]d\lambda \label{eq:rho}
\end{equation}
can be approximated by a Gaussian quadrature rule on $[0,1]$, where the quadrature points $\{\lambda_i\}$ are mapped to $\{\xi_i\}$ and $\{\eta_i\}$ by~\eqref{eq:lambda}.

\subsection{Spatial Discretization}
\begin{figure}[htb]
\begin{subfigure}{0.32\textwidth}
	\begin{tikzpicture} [x=1.5cm,y=1.25cm]
    \foreach \i in {0,...,2} {
        \draw [black, thick] (\i,0) -- (\i,2);
        \draw [black, thick] (0,\i) -- (2,\i) ;
    }
    \foreach \i in {1,...,2} {
        \draw [black, dashed] (\i-1/2,0) -- (\i-1/2,2) ;
        \draw [black, dashed] (0,\i-1/2) -- (2,\i-1/2);
    }
   \foreach \x in {1,...,2}{
      \foreach \y in {1,...,2}{
        \node[draw,circle,inner sep=2.5pt,fill,red] at (\x-1/2,\y-1/2) {};
      }
    }
    \foreach \x in {0,...,2}{
      \foreach \y in {1,...,2}{
        \node[draw,diamond,inner sep = 2.pt,fill,blue] at (\x,\y-1/2) {};
      }
    }
    \foreach \x in {1,...,2}{
      \foreach \y in {0,...,2}{
        \node[draw,diamond,inner sep = 2.pt,fill,blue] at (\x-1/2,\y) {};
      }
    }
        \foreach \x in {0,...,2}{
      \foreach \y in {0,...,2}{
        \node[draw,circle,inner sep=2.5pt,fill,red] at (\x,\y) {};
      }
    }
    \node [below] at (0,0) {$x_{i-1}$};
    \node [below] at (1,0) {$x_{i}$};
    \node [below] at (2,0) {$x_{i+1}$};
    \node [left] at (0,0) {$y_{j-1}$};
    \node [left] at (0,1) {$y_{j}$};
    \node [left] at (0,2) {$y_{j+1}$};
	\end{tikzpicture}
\caption{Staggered grids} \label{fig:CartesianGrid}
\end{subfigure}
\centering
\begin{subfigure}{0.32\textwidth}
	\begin{tikzpicture} [x=1.5cm,y=1.25cm]
 \foreach \i in {0,...,2} {
        \draw [thin,gray!70!white] (\i,0) -- (\i,2);
        \draw [thin,gray!70!white] (0,\i) -- (2,\i) ;
    }
    \foreach \i in {1,...,2} {
        \draw [thin,gray!70!white,dashed] (\i-1/2,0) -- (\i-1/2,2) ;
        \draw [thin,gray!70!white,dashed] (0,\i-1/2) -- (2,\i-1/2);
    }
    \foreach \i in {0,1} {
        \draw [black,thick] (0,\i+1/2) -- (2-1/2-\i,2);
        \draw [black,thick] (\i+1/2,0) -- (2,1.5-\i) ;
        \draw [black,thick] (\i+1/2,0) -- (0,\i+1/2);
        \draw [black,thick] (2,\i+1/2) -- (\i+1/2,2) ;
    }
   \foreach \x in {1,...,2}{
      \foreach \y in {1,...,2}{
        \node[draw,circle,inner sep=2.5pt,fill,red] at (\x-1/2,\y-1/2) {};
      }
    }
    \foreach \x in {0,...,2}{
      \foreach \y in {1,...,2}{
        \node[draw,diamond,inner sep = 2.pt,fill,blue] at (\x,\y-1/2) {};
      }
    }
    \foreach \x in {1,...,2}{
      \foreach \y in {0,...,2}{
        \node[draw,diamond,inner sep = 2.pt,fill,blue] at (\x-1/2,\y) {};
      }
    }
        \foreach \x in {0,...,2}{
      \foreach \y in {0,...,2}{
        \node[draw,circle,inner sep=2.5pt,fill,red] at (\x,\y) {};
      }
    }
    \node [below] at (0,0) {$x_{i-1}$};
    \node [below] at (1,0) {$x_{i}$};
    \node [below] at (2,0) {$x_{i+1}$};
    \node [left] at (0,0) {$y_{j-1}$};
    \node [left] at (0,1) {$y_{j}$};
    \node [left] at (0,2) {$y_{j+1}$};
\end{tikzpicture}
\caption{control volumes of 
\begin{tikzpicture}
\node[draw,circle,inner sep=2.5pt,fill,red] at (0,0) {};
\end{tikzpicture}
} \label{fig:controlVolume1}
\end{subfigure}
\begin{subfigure}{0.32\textwidth}
	\begin{tikzpicture} [x=1.5cm,y=1.25cm]
 \foreach \i in {0,...,2} {
        \draw [thin,gray!70!white] (\i,0) -- (\i,2);
        \draw [thin,gray!70!white] (0,\i) -- (2,\i) ;
    }
    \foreach \i in {1,...,2} {
        \draw [thin,gray!70!white,dashed] (\i-1/2,0) -- (\i-1/2,2) ;
        \draw [thin,gray!70!white,dashed] (0,\i-1/2) -- (2,\i-1/2);
    }
    \foreach \i in {0,1} {
        \draw [black,thick] (0,\i) -- (2-\i,2);
        \draw [black,thick] (1,0) -- (2,1) ;
        \draw [black,thick] (0,2-\i) -- (2-\i,0);
        \draw [black,thick] (1,2) -- (2,1) ;
    }
   \foreach \x in {1,...,2}{
      \foreach \y in {1,...,2}{
        \node[draw,circle,inner sep=2.5pt,fill,red] at (\x-1/2,\y-1/2) {};
      }
    }
    \foreach \x in {0,...,2}{
      \foreach \y in {1,...,2}{
        \node[draw,diamond,inner sep = 2.pt,fill,blue] at (\x,\y-1/2) {};
      }
    }
    \foreach \x in {1,...,2}{
      \foreach \y in {0,...,2}{
        \node[draw,diamond,inner sep = 2.pt,fill,blue] at (\x-1/2,\y) {};
      }
    }
        \foreach \x in {0,...,2}{
      \foreach \y in {0,...,2}{
        \node[draw,circle,inner sep=2.5pt,fill,red] at (\x,\y) {};
      }
    }
    \node [below] at (0,0) {$x_{i-1}$};
    \node [below] at (1,0) {$x_{i}$};
    \node [below] at (2,0) {$x_{i+1}$};
    \node [left] at (0,0) {$y_{j-1}$};
    \node [left] at (0,1) {$y_{j}$};
    \node [left] at (0,2) {$y_{j+1}$};
\end{tikzpicture}
\caption{control volumes of
\begin{tikzpicture}
\node[draw,diamond,inner sep = 2.pt,fill,blue] at (0,0) {};
\end{tikzpicture}
} \label{fig:controlVolume2}
\end{subfigure}
\caption{Staggered grids and control volumes:
 red circles (vertices and cell centers), $\pr{1},\; \pr{2},\; \rho,\; Q$; blue diamonds (face centers), $\pj{1},\;\pj{2}$.}
\label{fig:StaggeredGrids}
\end{figure}
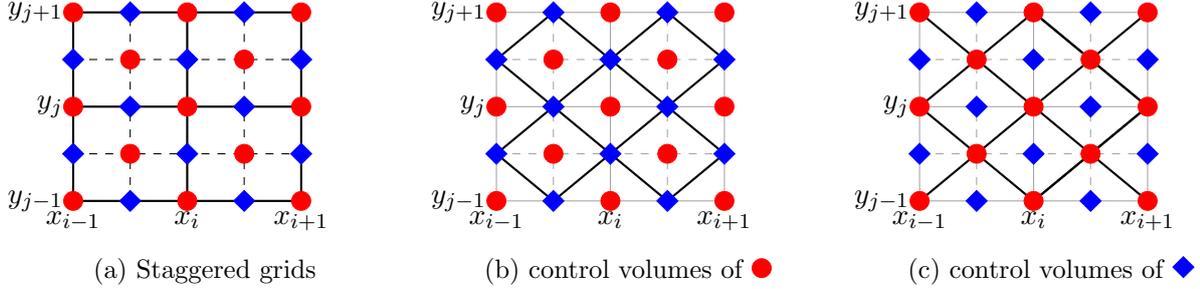

First, we describe the spatial discretization, while keeping the time $t$ continuous. We define a standard regular mesh with mesh size $\Delta x \times \Delta y$ and place the unknowns in the following way (see Figure~\ref{fig:CartesianGrid}):
\begin{itemize}
\item $\pr{1},\; \pr{2},\; \rho$, and $Q$ are located at the vertices $(i,j)$ and the cell centers $(i+\frac12,j+\frac12)$;
\item $\pj{1}$ and $\pj{2}$ are located at the face centers $(i+\frac12,j)$ and $(i,j+\frac12)$.
\end{itemize}
This choice enables us to approximate all spatial derivatives in the system~\eqref{eq:pde_system} by half-grid centered finite differences and yields a closed system of equations.

In detail, the semidiscretized equations are defined as follows. Let $i,j\in \mathbb{Z}$. Then, the parities $\pr{1}$ and $\pr{2}$ satisfy~\eqref{eq:r1} and~\eqref{eq:r2}, which on the vertices $(i,j)$ and the cell centers $(i+\frac12,j+\frac12)$ are given by
\begin{equation}
\begin{aligned}\label{eq:discret_r1}
\partial_t \pr{1}_{i,j} + \xi \frac{\pj{1}_{i+\frac 1 2,j}-\pj{1}_{i-\frac 1 2,j}}{\Delta x} - \eta \frac{\pj{1}_{i,j+\frac 1 2}-\pj{1}_{i,j-\frac 1 2}}{\Delta y} 
= -\frac{\sigma_s}{\varepsilon^2} (\pr{1}_{i,j}-\rho_{i,j}) - \sigma_a \pr{1}_{i,j} + Q_{i,j} \;,\\
\partial_t \pr{1}_{i+\frac 1 2,j+ \frac 1 2} + \xi \frac{\pj{1}_{i+1,j+\frac 1 2}- \pj{1}_{i,j+\frac 1 2}}{\Delta x} - \eta \frac{\pj{1}_{i+\frac 1 2,j+1}-\pj{1}_{i+\frac 1 2,j}}{\Delta y} 
=-\frac{\sigma_s}{\varepsilon^2} (\pr{1}_{i+\frac 1 2, j+\frac 1 2}-\rho_{i+\frac 1 2,j+\frac 1 2}) \\
-\sigma_a \pr{1}_{i+\frac 1 2,j+\frac 1 2} + Q_{i+\frac 1 2,j+\frac 1 2} \;,
\end{aligned}
\end{equation}
and
\begin{equation}
\begin{aligned}\label{eq:discret_r2}
\partial_t \pr{2}_{i,j} + \xi \frac{\pj{2}_{i+\frac 1 2,j}-\pj{2}_{i-\frac 1 2,j}}{\Delta x} + \eta \frac{\pj{2}_{i,j+\frac 1 2}-\pj{2}_{i,j-\frac 1 2}}{\Delta y} 
= -\frac{\sigma_s}{\varepsilon^2} (\pr{2}_{i,j}-\rho_{i,j}) - \sigma_a \pr{2}_{i,j} + Q_{i,j} \;,\\
\partial_t \pr{2}_{i+\frac 1 2,j+ \frac 1 2} + \xi \frac{\pj{2}_{i+1,j+\frac 1 2}- \pj{2}_{i,j+\frac 1 2}}{\Delta x} + \eta \frac{\pj{2}_{i+\frac 1 2,j+1}-\pj{2}_{i+\frac 1 2,j}}{\Delta y} 
=-\frac{\sigma_s}{\varepsilon^2} (\pr{2}_{i+\frac 1 2, j+\frac 1 2}-\rho_{i+\frac 1 2,j+\frac 1 2}) \\
-\sigma_a \pr{2}_{i+\frac 1 2,j+\frac 1 2} + Q_{i+\frac 1 2,j+\frac 1 2} \;.
\end{aligned}
\end{equation}

Similarly, the equations for the parities $\pj{1}$ and $\pj{2}$ on the face centers $(i+\frac 1 2,j)$ and $(i,j+\frac 1 2)$ are derived from the equations~\eqref{eq:j1} and~\eqref{eq:j2}:
\begin{equation}
\begin{aligned}\label{eq:discret_j1}
\partial_t \pj{1}_{i+\frac 1 2,j} 
+ \frac\xi{\varepsilon^2} \frac{\pr{1}_{i+1,j}- \pr{1}_{i,j}}{\Delta x}
-\frac\eta{\varepsilon^2} \frac{\pr{1}_{i+\frac 1 2,j+\frac 1 2}-\pr{1}_{i+\frac 1 2,j-\frac 1 2}}{\Delta y} 
&=  -\frac{\sigma_t}{\varepsilon^2}  \pj{1}_{i+\frac 1 2,j} \;,\\
\partial_t \pj{1}_{i,j+\frac 1 2} 
+ \frac\xi{\varepsilon^2} \frac{\pr{1}_{i+\frac 1 2,j+\frac 1 2}- \pr{1}_{i-\frac 1 2,j+\frac 1 2}}{\Delta x}
- \frac\eta{\varepsilon^2} \frac{\pr{1}_{i,j+1}-\pr{1}_{i,j}}{\Delta y}
&= -\frac{\sigma_t}{\varepsilon^2} \pj{1}_{i,j+\frac 1 2}\;,
\end{aligned}
\end{equation}
and
\begin{equation}
\begin{aligned}\label{eq:discret_j2}
\partial_t \pj{2}_{i+\frac 1 2,j} 
+ \frac\xi{\varepsilon^2} \frac{\pr{2}_{i+1,j}- \pr{2}_{i,j}}{\Delta x}
+ \frac\eta{\varepsilon^2} \frac{\pr{2}_{i+\frac 1 2,j+\frac 1 2}-\pr{2}_{i+\frac 1 2,j-\frac 1 2}}{\Delta y} 
&=  -\frac{\sigma_t}{\varepsilon^2}  \pj{2}_{i+\frac 1 2,j} \;,\\
\partial_t \pj{2}_{i,j+\frac 1 2} 
+ \frac\xi{\varepsilon^2} \frac{\pr{2}_{i+\frac 1 2,j+\frac 1 2}- \pr{2}_{i-\frac 1 2,j+\frac 1 2}}{\Delta x}
+ \frac\eta{\varepsilon^2} \frac{\pr{1}_{i,j+1}-\pr{1}_{i,j}}{\Delta y}
&= -\frac{\sigma_t}{\varepsilon^2} \pj{2}_{i,j+\frac 1 2}\;.
\end{aligned}
\end{equation}

\begin{remark}
The method can be interpreted as a finite volume method and is therefore conservative. To see this, note that all the unknowns are given on two regular grids: vertices and cell centers or vertical cell centers and horizontal cell centers. Merging the corresponding grids, we obtain two staggered nonregular grids as shown in Figure~\ref{fig:CartesianGrid}. These can also be interpreted as control volumes of a finite volume method. The control volumes $\Diamond_{i,j}$ are defined as shown in Figures~\ref{fig:controlVolume1}, \ref{fig:controlVolume2}, where $\Diamond_{i,j}$ is the diamond around the point $(x_i,y_j)$.

We define the volume averages corresponding to the control volumes by, e.g.,
\begin{align}
\pr{1}_{i,j} := \frac{1}{|\Diamond_{i,j}|} \int_{\Diamond_{i,j}} \pr{1} d(x,y)\quad  \text{with}\quad  |\Diamond_{i,j}| = \frac12 \Delta x\Delta y
\end{align}
and integrate the system~\eqref{eq:pde_system} over the control volumes. For instance,~\eqref{eq:r1} integrated over the volume $\Diamond_{i,j}$ is given by
\begin{align}
\partial_t \pr{1}_{i,j} +\frac{1}{|\Diamond_{i,j} | } \int_{\Diamond_{i,j}} \left(\xi \partial_x \pj{1} + \eta \partial_y \pj{1}\right) d(x,y) = -\frac{\sigma_s}{\varepsilon^2} (\pr{1}_{i,j}- \rho_{i,j}) +Q_{i,j}\;.
\end{align}
The integral can be simplified using  Gauss's theorem
\begin{equation}
\begin{aligned}
\int_{\Diamond_{i,j}} \left(\xi \partial_x \pj{1} + \eta \partial_y \pj{1}\right) d(x,y)
 &= \int_{\Diamond_{i,j}} \left(\nabla\cdot \begin{bmatrix} \xi \\ \eta\end{bmatrix} \pj{1} \right)d(x,y) \\
 &= \int_{\partial\Diamond_{i,j}} \left(\begin{bmatrix} \xi \\ \eta\end{bmatrix} \pj{1} \cdot n \right) d(x,y)\;,
\end{aligned}
\end{equation}
where $n$ is the outer normal vector. It remains to compute the integrals over the four edges of the diamond. These terms are approximated by the trapezoidal rule, e.g., the upper right part is then given by
\begin{multline}
\int_{(x_{i+\frac12},y_j)}^{(x_i,y_{j+\frac12})} \left(\begin{bmatrix} \xi \\ \eta\end{bmatrix} \pj{1} \cdot n \right) d(x,y)
= \int_{(x_{i+\frac12},y_j)}^{(x_i,y_{j+\frac12})}\left(\begin{bmatrix} \xi \\ \eta\end{bmatrix} \pj{1} \cdot \begin{bmatrix} \tfrac{\Delta y}{\sqrt{\Delta x^2+ \Delta y^2}} \\ \tfrac{\Delta x }{\sqrt{\Delta x^2+ \Delta y^2}} \end{bmatrix}\right) d(x,y)\\
 \approx \frac12\left( \xi \Delta y 
 \left(\pj{1}_{i,j+\frac12} +\pj{1}_{i+\frac12,j}\right)  +  \eta \Delta x \left(\pj{1}_{i,j+\frac12} +\pj{1}_{i+\frac12,j} \right)\right)\;.
\end{multline}
In total, we obtain the approximation for the integral
\begin{align}
\frac{2}{\Delta x\Delta y } \int_{\partial\Diamond_{i,j}} \left(\begin{bmatrix} \xi \\ \eta\end{bmatrix} \pj{1} \cdot n \right) d(x,y)
\approx \frac{\xi}{\Delta x}\left( \pj{1}_{i+\frac12,j} - \pj{1}_{i-\frac12,j} \right) +  \frac{\eta}{\Delta y}\left( \pj{1}_{i,j+\frac12} - \pj{1}_{i,j-\frac12}  \right)\;,
\end{align}
and therefore the same semidiscrete equation as before~\eqref{eq:discret_r1}
\begin{align}
\partial_t \pr{1}_{i,j} + \frac{\xi}{\Delta x}\left( \pj{1}_{i+\frac12,j} - \pj{1}_{i-\frac12,j} \right) +  \frac{\eta}{\Delta y}\left( \pj{1}_{i,j+\frac12} - \pj{1}_{i,j-\frac12}  \right) = -\frac{\sigma_s}{\varepsilon^2} (\pr{1}_{i,j}- \rho_{i,j}) +Q_{i,j}\;.
\end{align}
In the same way, we obtain the semidiscretized equations~\eqref{eq:discret_r2}-\eqref{eq:discret_j2}.

In the case of a Cartesian mesh $\Delta x= \Delta y$, the finite volumes are rotated by 45 degrees with respect to the axes, which seems to be odd. However, it seems that the choice of grid points and therefore the volumes are unique in the strategy we have adopted. In fact, it seems quite natural when we consider the fact that the quantity $\pr{2}$ describes the number of particles in $\xi,\eta>0$ and $\xi,\eta<0$. Thus a flow into the diagonal direction makes sense, and therefore interfaces at 45 degrees to the axes. This will be investigated further in future work.
\end{remark}

\begin{remark}
Note that \eqref{eq:definition-parities} can formally be inverted to obtain the density of the transport equation
\begin{align}
f(\xi,\eta) = \begin{cases} \pr{1}(|\xi|,|\eta|)+\varepsilon\;\text{sign}(\xi)\pj{1}(|\xi|,|\eta|)& \text{for}\; \xi\eta<0\;,\\
\pr{2}(|\xi|,|\eta|)+\varepsilon\;\text{sign}(\xi)\pj{2}(|\xi|,|\eta|) &\text{for}\;\xi\eta\geq0\;.
\end{cases}
\end{align}
However, in the numerical scheme the parities $\pr{1}$, $\pr{2}$ and $\pj{1}$, $\pj{2}$ are not given on the same spatial grid and need to be interpolated. This is similar to other approaches which are based on parity decompositions.
\end{remark}

\subsection{Time Discretization}
For simplicity, we consider again semidiscretized equations. This time, we keep the spatial variables $x$ and $y$ continuous and apply the time discretization technique from~\cite{Jin-Pareschi-Toscani-2000}. The idea is to introduce a relaxation parameter $\phi = \phi(\varepsilon)$, such that we obtain a linear hyperbolic system with stiff relaxation. Then the linear hyperbolic part, which is nonstiff, can be separated from the stiff relaxation step.

First, we rewrite the system of~\eqref{eq:pde_system} as the diffusive relaxation system
\begin{equation}
\begin{aligned}
\partial_t \pr{1} +\xi \partial_x \pj{1} -\eta \partial_y \pj{1} &=-\frac{\sigma_s}{\varepsilon^2} (\pr{1} - \rho)-\sigma_a \pr{1}+Q \;,\\
\partial_t \pr{2} +\xi \partial_x \pj{2} +\eta \partial_y \pj{2} &=-\frac{\sigma_s}{\varepsilon^2} (\pr{2} - \rho)-\sigma_a \pr{2}+Q \;,\\
\partial_t \pj{1} + \phi\xi \partial_x \pr{1} -\phi\eta \partial_y \pr{1} 
&=-\frac{1}{\varepsilon^2}[\sigma_s \pj{1}+(1-\varepsilon^2\phi)\xi\partial_x \pr{1} - (1-\varepsilon^2 \phi) \eta \partial_y \pr{1}] \;,\\
\partial_t \pj{2} +\phi\xi\partial_x \pr{2} +\phi \eta \partial_y \pr{2}
&=-\frac{1}{\varepsilon^2}[\sigma_s \pj{2}+ (1-\varepsilon^2 \phi)\xi \partial_x \pr{2} + (1-\varepsilon^2 \phi) \eta \partial_y \pr{2}] 
\end{aligned}
\end{equation}
with $0\leq \phi\leq 1/\varepsilon^2$. The condition $\phi\geq 0$ is necessary for the hyperbolicity of the left-hand side, whereas the condition $\phi\leq 1/\varepsilon^2$ ensures that the bracketed term on the right-hand side has a well-defined limit for $\varepsilon\to 0$. In Remark \ref{rem:phi} below we comment on the role of $\phi$ and how this choice differs from Klar's scheme~\cite{Klar-1998}.

Thus, we split the equation into two parts, the transport step,
\begin{equation}\label{eq:transport_step}
\begin{aligned}
\partial_t \pr{1} + \xi \partial_x \pj{1} -\eta \partial_y \pj{1} &= -\sigma_a \pr{1}+Q \;,\\
\partial_t \pr{2} + \xi \partial_x \pj{2} +\eta \partial_y \pj{2} &= -\sigma_a \pr{2}+Q \;,\\
\partial_t \pj{1} + \phi\xi \partial_x \pr{1} -\phi\eta \partial_y \pr{1} &=-\sigma_a \pj{1}\;,\\
\partial_t \pj{2} +\phi\xi\partial_x \pr{2} +\phi \eta \partial_y \pr{2} &= -\sigma_a \pj{2} \;,\\
\end{aligned}
\end{equation}
and the relaxation step,
\begin{equation}\label{eq:relaxation_step}
\begin{aligned}
\partial_t \pr{1}& = -\frac{\sigma_s}{\varepsilon^2} (\pr{1} - \rho) \;,\\
\partial_t \pr{2} &= -\frac{\sigma_s}{\varepsilon^2} (\pr{2} - \rho) \;,\\
\partial_t \pj{1} &= -\frac 1 {\varepsilon^2} [\sigma_s \pj{1}+ (1-\varepsilon^2 \phi)\xi \partial_x \pr{1} - (1-\varepsilon^2 \phi) \eta \partial_y \pr{1}] \;,\\
\partial_t \pj{2}  &= -\frac 1 {\varepsilon^2} [\sigma_s \pj{2}+ (1-\varepsilon^2 \phi)\xi \partial_x \pr{2} + (1-\varepsilon^2 \phi) \eta \partial_y \pr{2}] \;.
\end{aligned}
\end{equation}
Finally, we apply the explicit Euler method to the first step and the implicit Euler method to the second step. Note that the implicit Euler method can be implemented explicitly, since $\rho$ is preserved in the second step (which can be seen by adding the first two equations).

The fully discrete scheme is just splitting (\ref{eq:discret_r1})-(\ref{eq:discret_j2}) into the two steps (\ref{eq:transport_step})-(\ref{eq:relaxation_step}).

\begin{remark}
\label{rem:phi}
Klar~\cite{Klar-1998} developed a similar decomposition, which corresponds to $\phi=0$ in our framework. As mentioned in \cite{Jin-Pareschi-Toscani-2000}, there are two major differences. First, the resulting system is only weakly hyperbolic and therefore well-posedness is an issue. Second, computations are still performed on the hole velocity space. Using the symmetries of the parities, the computations can be performed on the first quadrant and the computational cost can be reduced.
\end{remark}

\begin{remark}
To generalize the problem from a two-dimensional ``flatland'' to a full two-dimensional problem, the domain of the angular variable changes from the unit circle to the unit disc. In order to do this change, the Gauss integration rule needs to be substituted by a two-dimensional integration rule on the unit circle.
\end{remark}

\section{The AP property}
\label{sec:AP-property}
In this section, we analyze the AP property of the above scheme in two steps. First, we derive the discrete asymptotic limit. Second, we analyze stability.

\subsection{The Diffusion limit}
\label{sec:diff-limit}
In the same way as above, we consider the spatial and the time discretization separately. The limit, as $\varepsilon \to 0$,  of the time discretization is derived in~\cite{Jin-Pareschi-Toscani-2000}. Hence, it remains to investigate the discrete limit of the spatial discretization. To this end, we consider the diffusive limit $\varepsilon\to0$ of the semidiscretized equations~\eqref{eq:discret_r1}-\eqref{eq:discret_j2}.

First, the limit of~\eqref{eq:discret_r1} and~\eqref{eq:discret_j1} for the parities $\pr{1}$ and $\pj{1}$ is given by
\begin{equation}
\begin{aligned}
\pr{1}_{i,j} &= \rho_{i,j}\;, \\
\pr{1}_{i+\frac 1 2,j+\frac 1 2} &= \rho_{i+\frac 1 2,j+\frac 1 2}\;, \\
\pj{1}_{i+\frac 1 2,j}&= - \frac\xi{\sigma_t} \frac{\pr{1}_{i+1,j}- \pr{1}_{i,j}}{\Delta x} +\frac\eta{\sigma_s} \frac{\pr{1}_{i+\frac 1 2,j+\frac 1 2}-\pr{1}_{i+\frac 1 2,j-\frac 1 2}}{\Delta y} \;,\\
\pj{1}_{i,j+\frac 1 2} &= -\frac\xi{\sigma_t}  \frac{\pr{1}_{i+\frac 1 2,j+\frac 1 2}- \pr{1}_{i-\frac 1 2,j+\frac 1 2}}{\Delta x}+ \frac\eta{\sigma_s} \frac{\pr{1}_{i,j+1}-\pr{1}_{i,j}}{\Delta y}\;.
\end{aligned}
\end{equation}
Inserting these equations into~\eqref{eq:discret_r1} yields
\begin{equation}
\begin{aligned}\label{eq:rho_r1}
 \partial_t \rho_{i,j} 
&- \frac{\xi^2}{\sigma_t} \frac{\rho_{i+1,j}-2\rho_{i,j}+\rho_{i-1,j}}{(\Delta x)^2} \\
&+ \frac{2\xi \eta}{\sigma_t} \frac{\rho_{i+\frac 1 2, j+\frac 1 2} - \rho_{i+\frac 1 2,j-\frac 1 2}-\rho_{i-\frac 1 2,j+\frac 1 2}+\rho_{i-\frac 1 2,j-\frac 1 2}}{\Delta x \Delta y} \\
 &-\frac{\eta^2}{\sigma_t} \frac{\rho_{i,j+1}-2\rho_{i,j}+\rho_{i,j-1}}{(\Delta y)^2}
=-\sigma_a \rho_{i,j} + Q_{i,j} \;,\\
 \partial_t \rho_{i+ \frac 1 2,j+\frac 1 2} 
&- \frac{\xi^2}{\sigma_t} \frac{\rho_{i+ \frac 3 2, j+\frac 1 2} - 2 \rho_{i+\frac 1 2,j+\frac 1 2} + \rho_{i-\frac 1 2,j+\frac 1 2}}{(\Delta x)^2} \\
&+ \frac{2\xi \eta}{\sigma_t} \frac { \rho_{i+1,j+1}-\rho_{i+1,j}-\rho_{i,j+1}+\rho_{i,j}}{\Delta x\Delta y } \\
&- \frac{\eta^2}{\sigma_t} \frac{\rho_{i+\frac 1 2,j+\frac 3 2}- 2 \rho_{i+\frac 1 2,\frac 1 2}+\rho_{i+\frac 1 2,j-\frac 1 2} }{(\Delta y)^2}
= -\sigma_a \rho_{i+\frac 1 2, j+ \frac 1 2} + Q_{i+\frac 1 2, j+ \frac 1 2}\;.
\end{aligned}
\end{equation}
Treating~\eqref{eq:discret_r2} and~\eqref{eq:discret_j2} in the same way as above, we additionally obtain the following differential equations for $\rho$:
\begin{equation}
\begin{aligned}\label{eq:rho_r2}
\partial_t \rho_{i,j} 
 &- \frac{\xi^2}{\sigma_s} \frac{\rho_{i+1,j}-2\rho_{i,j}+\rho_{i-1,j}}{(\Delta x)^2} \\
& - \frac{2\xi \eta}{\sigma_s} \frac{\rho_{i+\frac 1 2, j+\frac 1 2} - \rho_{i+\frac 1 2,j-\frac 1 2}-\rho_{i-\frac 1 2,j+\frac 1 2}+\rho_{i-\frac 1 2,j-\frac 1 2}}{\Delta x \Delta y} \\
 & -\frac{\eta^2}{\sigma_s} \frac{\rho_{i,j+1}-2\rho_{i,j}+\rho_{i,j-1}}{(\Delta y)^2}
 = -\sigma_a \rho_{i,j} + Q_{i,j} \;,\\
\partial_t \rho_{i+ \frac 1 2,j+\frac 1 2} 
& - \frac{\xi^2}{\sigma_s} \frac{\rho_{i+ \frac 3 2, j+\frac 1 2} - 2 \rho_{i+\frac 1 2,j+\frac 1 2} + \rho_{i-\frac 1 2,j+\frac 1 2}}{(\Delta x)^2} \\
& - \frac{2\xi \eta}{\sigma_s} \frac { \rho_{i+1,j+1}-\rho_{i+1,j}-\rho_{i,j+1}+\rho_{i,j}}{\Delta x\Delta y } \\
 &- \frac{\eta^2}{\sigma_s} \frac{\rho_{i+\frac 1 2,j+\frac 3 2}- 2 \rho_{i+\frac 1 2,\frac 1 2}+\rho_{i+\frac 1 2,j-\frac 1 2} }{(\Delta y)^2}
 = -\sigma_a \rho_{i+\frac 1 2, j+ \frac 1 2} + Q_{i+\frac 1 2, j+ \frac 1 2} \;.
\end{aligned}
\end{equation}
Adding up the equations, the middle terms cancel and we obtain
\begin{equation}
\begin{aligned}
\begin{aligned}
\partial_t \rho_{i,j}
 - \frac{\xi^2}{\sigma_s} \frac{\rho_{i+1,j}-2\rho_{i,j}+\rho_{i-1,j}}{(\Delta x)^2} 
 -\frac{\eta^2}{\sigma_s} \frac{\rho_{i,j+1}-2\rho_{i,j}+\rho_{i,j-1}}{(\Delta y)^2}
 = -\sigma_a \rho_{i,j} + Q_{i,j} \;,
\end{aligned}\\
\begin{aligned}
\partial_t \rho_{i+ \frac 1 2,j+\frac 1 2} 
 &- \frac{\xi^2}{\sigma_s} \frac{\rho_{i+ \frac 3 2, j+\frac 1 2} - 2 \rho_{i+\frac 1 2,j+\frac 1 2} + \rho_{i-\frac 1 2,j+\frac 1 2}}{(\Delta x)^2}\\
 &- \frac{\eta^2}{\sigma_s} \frac{\rho_{i+\frac 1 2,j+\frac 3 2}- 2 \rho_{i+\frac 1 2,\frac 1 2}+\rho_{i+\frac 1 2,j-\frac 1 2} }{(\Delta y)^2}
 = -\sigma_a \rho_{i+\frac 1 2, j+ \frac 1 2} + Q_{i+\frac 1 2, j+ \frac 1 2} \;.
\end{aligned}
\end{aligned}
\end{equation}
Integrating over $\xi^2+\eta^2=1$ yields the semidiscretized diffusion equations on the vertices and the cell centers. Note that integrating~(\ref{eq:rho_r1}) or~(\ref{eq:rho_r2}) over $\xi^2+\eta^2=1$, the middle terms cancel as well and we get the same result.

As expected, the spatial discretization with staggered grids leads to a compact five-point stencil for the diffusion equation~\eqref{eq:diffusion}. Together with the results~\cite{Jin-Pareschi-Toscani-2000} on the limit of the time discretization~\eqref{eq:transport_step} --~\eqref{eq:relaxation_step}, this also shows that the formal limit of our scheme coincides with the diffusion equation.

\subsection{Stability}
\label{sec:stability}
We limit our discussion to the one-dimensional case (see Remark~\ref{rem:stability2D} for the two-dimensional case) and show uniform stability with $\varepsilon$ using the von Neumann analysis~\cite{Trefethen-1996,Seibold-Frank-2012,Lemou-Mieussens-2008}. 

In the following, we consider the transport equation in slab geometry and assume that the cross section $\sigma_t=\sigma_s+\varepsilon^2\sigma_a>0$ is independent of $x\in\R$ (see Remark~\ref{rem:x-dependence} for space-dependent scattering). Further, we consider a source-free two velocity model $v\in\{\pm1\}$. Then, the even and odd parities
\begin{equation}
\pri(t,x)=\frac1 2 [f( t,x,1)+f(t,x,-1)]  \quad\text{and}\quad \pji(t,x)= \frac 1 {2\varepsilon} [f(t,x,1)-f(t,x,-1)]
\end{equation}
 fulfill 
\begin{equation}
\begin{aligned}
\partial_t \pri + \partial_x \pji &=-\sigma_a \pri \;,\\
\partial_t \pji + \tfrac 1 {\varepsilon^2} \partial_x \pri &= -\tfrac 1 {\varepsilon^2} \sigma_s \pji  - \sigma_a \pji  \;,
\end{aligned}
\end{equation}
and the numerical scheme has the following update rule: For $k=0,1,2,\ldots$
\begin{equation}
\begin{aligned}
\pri^{k+\frac 1 2} &= \pri^k - \Delta t( D_x \pji^k + \sigma_a \pri^k) \;,\\
\pji^{k+\frac 1 2} &= \pji ^k - \Delta t( \phi D_x \pri^k+ \sigma_a \pji^k) \;,\\ 
\pri^{k+1} &= \pri^{k+\tfrac 1 2} \;,\\
\pji^{k+1}  &= \tfrac{\varepsilon^2}{\varepsilon^2+\sigma_s\Delta t} \pji^{k+\tfrac 1 2 } - \tfrac{\Delta t}{\varepsilon^2+\sigma_s\Delta t} (1-\varepsilon^2 \phi) D_x \pri^{k+1} \;,
\end{aligned}
\end{equation}
where $D_x$ denotes the half-grid centered finite difference of the spatial derivative. We place $\pri$ on the half grid points $(m+\tfrac 1 2) \Delta x$ and $\pji$ on the full grid points $m\Delta x$. 
For this scheme, we do not expect positivity, but we seek a uniform CFL condition.

For a von Neumann analysis of the scheme, we expand the parities in Fourier series: 
\begin{equation}
\pri(x,t)= \sum_{\ell=-\infty}^\infty a_\ell(t) e^{i\ell x} \quad \text{and} \quad \pji(x,t)= \sum_{\ell-\infty}^\infty b_\ell(t) e^{i\ell x}\;.
\end{equation}
As no mixing between the Fourier modes occurs during the update of the solution, it is sufficient to consider the evolution of
\begin{equation}
 \pri(x,t)= a_\ell(t)e^{i\ell x}  \quad \text{and} \quad \pji(x,t)=b_\ell(t)e^{i\ell x}
\end{equation}
for some $\ell$ and to determine the growth factor matrix of the Fourier coefficients. First, we note that the staggered grid derivatives can be rewritten as
\begin{equation}
\begin{aligned}
(D_x \pri)\left(h\left(m+\tfrac 1 2\right),t\right) &= a_\ell(t)\frac{e^{i \ell h(m+1)}-e^{i \ell hm}} h =2\tfrac i h \sin\left(\tfrac{\ell h}2\right) e^{i \ell h(m+\tfrac 1 2 )} a_\ell(t) \;,\\
(D_x \pji) (hm,t) &= b_\ell(t) \frac {e^{i\ell h(m+\tfrac 1 2)}- e^{i\ell h(m-\tfrac 1 2)}} h  = 2\tfrac i h \sin\left(\tfrac{\ell h}2\right) e^{i\ell hm} b_\ell(t) \;,
\end{aligned}
\end{equation}
with $h:=\Delta x$. To shorten the notation, we define $d_\ell:= \tfrac {2i}h \sin\left(\tfrac{\ell h}{2}\right)$. Then, the first update step of the Fourier coefficients is given by
\begin{equation}
\begin{aligned}
\begin{bmatrix}a_k \\ b_k\end{bmatrix}
(t+\Delta t) = 
\underbrace{
\begin{bmatrix}
1-\sigma_a\Delta t &-\Delta t d_l \\ - \Delta t \phi d_\ell &1-\sigma_a\Delta t
\end{bmatrix}}_{=:G_1}
\begin{bmatrix}a_k \\ b_k\end{bmatrix} (t)
\end{aligned}
\end{equation}
and the second step is given by
\begin{equation}
\begin{bmatrix}a_k \\ b_k\end{bmatrix}
(t+\Delta t) = 
\underbrace{\begin{bmatrix}
1 &0\\ -\tfrac{\Delta t}{\varepsilon^2+\sigma_s\Delta t}(1-\varepsilon^2\phi) d_\ell &\tfrac{\varepsilon^2}{\varepsilon^2+\sigma_s\Delta t}
\end{bmatrix}}_{=:G_2}
\begin{bmatrix}a_k \\ b_k\end{bmatrix} (t) \;.
\end{equation}
Thus, the growth factor matrix is $G:=G_2\cdot G_1$
\begin{equation}
G=
\begin{bmatrix}
1-\sigma_a\Delta t & -\Delta t d_\ell \\
 -\tfrac{d_\ell\Delta t}{\varepsilon^2+\sigma_s\Delta t}(\sigma_a\Delta t (1-\varepsilon^2\phi)+1) &
\tfrac{1}{\varepsilon^2+\sigma_s\Delta t}(\varepsilon^2(1-\sigma_a\Delta t) + d_\ell^2\Delta t^2 (1-\varepsilon^2\phi))
\end{bmatrix} \;.
\end{equation}
For stability, the eigenvalues of the matrix $G$ are of main interest. They can be written as
\begin{equation}
\lambda_{1,2} = g\pm \sqrt{g^2-\det(G)}
\end{equation}
with $g$ being the half trace and $\det(G)$ being the determinant of $G$:
\begin{equation}
\label{eq:half-trace_determinant}
\begin{aligned}
g&= \tfrac 1 2 \tfrac{1}{\varepsilon^2+\sigma_s\Delta t} \left(  \Delta t^2 d_\ell^2(1-\varepsilon^2\phi) + (1-\sigma_a \Delta t)(2\varepsilon^2+\sigma_s\Delta t)\right)\quad \text{and}\\
\det(G)&= \tfrac{\varepsilon^2}{\varepsilon^2+\sigma_s\Delta t}((1-\sigma_a \Delta t)^2 - \phi d_\ell^2\Delta t^2) \;.
\end{aligned}
\end{equation}

\begin{proposition}
\label{prop:main_stability_result}
Let the time step $\Delta t$ and the relaxation parameter $\phi$ satisfy
\begin{equation}\label{eq:condition_CFL}
\Delta t \leq\min\left\{\tfrac{1}{\sigma_a}, \max\{\tfrac{1}{2}\varepsilon h,\; \tfrac 1 4 h^2\sigma_t\}\right\}
\end{equation}
and
\begin{equation}\label{eq:condition_phi}
0\leq \phi \leq \begin{cases} \frac{h\sigma_t}{2\varepsilon^3}, & h \sigma_t\leq 2\varepsilon\;,\\
\frac{1}{\varepsilon^2}, &\text{otherwise.}\end{cases}
\end{equation}
Then, the numerical scheme is $L^2$-stable.
\end{proposition}

\begin{remark}
Note that the three terms in the time step restriction~\eqref{eq:condition_CFL} can be interpreted separately. First, the $\frac{1}{2}\varepsilon h $ term comes from the advection operator. Second, the $\frac14 h^2 \sigma_t$ term corresponds to the Courant limit for explicit diffusion. 
Third, the $\frac{1}{\sigma_a}$ term is a result of the explicit treatment of the corresponding relaxation term. In this case, we assume that this term is small, so that $\Delta t \leq \frac{1}{\sigma_a}$ is not the restrictive term. Otherwise the term could be treated implicitly, which would remove the restriction.
\end{remark}

\begin{remark}
Note that the above restriction on $\phi$ is stricter than the one suggested in~\cite{Jin-Pareschi-Toscani-2000} $0\leq \phi\leq \frac 1 {\varepsilon^2}$. Moreover, the condition $h \sigma_t\leq 2\varepsilon$ is satisfied, if and only if the hyperbolic condition $ \Delta t\leq \max\{\tfrac{1}{2}\varepsilon h,\; \tfrac 1 4 h^2\sigma_t\} = \tfrac{1}{2}\varepsilon h$ holds. This means, there are the following two cases:
\begin{equation}
\begin{array}{llll}
h \sigma_t\leq 2\varepsilon: & \Delta t \leq \min\left\{\tfrac{1}{\sigma_a}, \tfrac{1}{2}\varepsilon h\right\}
&\text{and}& 0\leq \phi \leq \tfrac{h\sigma_t}{2\varepsilon^3} \;,\\
h \sigma_t> 2\varepsilon: & \Delta t \leq \min\left\{\tfrac{1}{\sigma_a}, \tfrac 1 4 h^2\sigma_t\right\}
&\text{and}&0\leq \phi \leq \tfrac{1}{\varepsilon^2} \;.
\end{array}
\end{equation}
In addition, as $\varepsilon\to 0$ the time step restriction becomes $\Delta t \lesssim \min\{\tfrac{1}{\sigma_a},\tfrac 1 4 h^2\sigma_t\}$, which does not vanish.
\end{remark}

In the proof of the proposition, we use the von Neumann analysis. A complete overview of these stability conditions can be found in the lecture notes by Trefethen~\cite{Trefethen-1996}.

\begin{proof}
Stability follows from the von Neumann condition if we can show $|\lambda_{1,2}|\leq 1$ for $\lambda_1\neq\lambda_2$ and $|\lambda_{1,2}|<1$ for $\lambda_1=\lambda_2$. To show these inequalities, we consider three different cases: two complex eigenvalues, two real eigenvalues, and one eigenvalue. Since $g$ and $\det(G)$ are real-valued~\eqref{eq:half-trace_determinant}, the cases are equivalent to: $g^2<\det(G)$, $g^2>\det(G)$, and $g^2=\det(G)$.

\textbf{Case} $g^2<\det(G)$ (two complex eigenvalues):\quad
If the eigenvalues $\lambda_{1,2}$ are complex, their real part is $g$ and their imaginary part is $\pm\sqrt{\det(G)-g^2}$. Thus, the stability condition $|\lambda_{1,2}|^2\leq1$ is satisfied if $\det(G) \leq 1$. For the determinant we have the following estimate
\begin{equation}
\label{eq:estimate_det}
\det(G) = \tfrac{\varepsilon^2}{\varepsilon^2+\sigma_s\Delta t}((1-\sigma_a \Delta t)^2 - \phi d_\ell^2\Delta t^2)
 \leq  \tfrac{\varepsilon^2}{\varepsilon^2+\sigma_s\Delta t}(1-\sigma_a \Delta t+\phi\tfrac{4 \Delta t^2}{h^2}) \;,
\end{equation}
where we used that $-d_\ell^2 =\tfrac{4}{h^2}\sin^2(\tfrac{\ell h}{2})\leq \tfrac 4{h^2} $ and the CFL condition $\Delta t\leq \tfrac{1}{\sigma_a}$. It remains to show that the last term of~\eqref{eq:estimate_det} is bounded by $1$. This is equivalent to
\begin{equation}
 \varepsilon^2\phi\tfrac{4 \Delta t}{h^2} \leq \sigma_s+\varepsilon^2 \sigma_a =\sigma_t \;,
\end{equation}
which in turn is satisfied under the condition $\Delta t \leq \max\{\tfrac{1}{2}\varepsilon h,\; \tfrac 1 4 h^2\sigma_t\}$ and the assumption~\eqref{eq:condition_phi}. This is one of the reasons for the choice of the upper bound of $\phi$ in the assumption~\eqref{eq:condition_phi}.

\textbf{Case} $g^2>\det(G)$ (two real eigenvalues): \quad The determinant of $G$ is always positive and therefore the eigenvalues are either both positive or both negative, and their sign changes with the sign of $g$. Thus, it is sufficient to show $\lambda_1\leq 1$ if $g\geq0$ and $\lambda_2\geq -1$ if $g<0$.
In particular, one can show that this is equivalent to
\begin{equation}
\det(G)+1\mp 2g\geq0 \;.
\end{equation}
The first inequality is generic
\begin{equation}
\det(G)+1-2g = \tfrac{\Delta t^2}{\varepsilon^2+\sigma_s\Delta t}(\sigma_a^2\varepsilon^2+\sigma_s\sigma_a-d_\ell)\geq 0 \;,
\end{equation}
since $d_\ell^2\leq0$. Whereas, the second inequality requires the CFL condition~\eqref{eq:condition_CFL}. More precisely, under the condition $0\leq \phi\varepsilon^2\leq 1$ and $\Delta t \leq \tfrac{1}{\sigma_a}$, we obtain
\begin{multline}
\det(G)+1+2g > 1+2g
= 1 + \tfrac{1}{\varepsilon^2+\sigma_s \Delta t}\left( \Delta t^2d_\ell^2(1-\varepsilon^2 \phi)+ (1-\sigma_a \Delta t)(2\varepsilon^2+\sigma_s \Delta t)\right)\\
\geq \tfrac{1}{\varepsilon^2+\sigma_s\Delta t}(\varepsilon^2  - \tfrac{4\Delta t^2}{h^2}+ \sigma_s\Delta t) \;.
\end{multline}
On the one hand, this is obviously nonnegative under the condition $\Delta t \leq \tfrac 1 2 \varepsilon h$. On the other hand, the second term can be rewritten as
\begin{equation}
\varepsilon^2  - \tfrac{4\Delta t^2}{h^2}+ \sigma_s\Delta t =\varepsilon^2(1-\sigma_a\Delta t) +\Delta t( \sigma_t - \tfrac{4\Delta t}{h^2}) \;,
\end{equation}
which is nonnegative under the condition $\Delta t\leq \tfrac{1}{\sigma_a}$ and $\Delta t \leq \tfrac 1 4 h^2 \sigma_t$. Together, this yields the desired inequality $\det(G)+1+2g > 0$.

\textbf{Case} $g^2=\det(G)$ (one eigenvalue): \quad The eigenvalue of $G$ is $\lambda_1=\lambda_2=g$. Thus, we need to show $|g|<1$. But as $\det(G)+1+2g>0$ and $\det(G)\leq1$ (see above cases) already imply $g>-1$, it remains to show $g<1$. Since $\sigma_t= \sigma_s+\varepsilon^2\sigma_a> 0$, at least one of the terms $\sigma_a\Delta t$, $\sigma_s\Delta t$ is positive and we obtain
\begin{equation}
\begin{aligned}
g &\leq  \tfrac 1 2 \tfrac{1}{\varepsilon^2+\sigma_s\Delta t} \left( (1-\sigma_a \Delta t)(2\varepsilon^2+\sigma_s\Delta t)\right)\\
 &<  \tfrac 1 2 \tfrac{1}{\varepsilon^2+\sigma_s\Delta t} \left( (1-\sigma_a \Delta t +\sigma_a\Delta t)(2\varepsilon^2+\sigma_s\Delta t+\sigma_s\Delta t)\right) =1 \;.
 \end{aligned}
\end{equation}
\end{proof}

\begin{remark}
If there is neither scattering nor absorption and the conditions~\eqref{eq:condition_CFL} and~\eqref{eq:condition_phi} hold, then the relaxation parameter satisfies $\phi=0$ and the hyperbolic condition is always satisfied. Further, the determinant, $\det(G) = 1$, and the half trace, $ g = 1+\frac{\Delta t^2 d_\ell^2}{2\varepsilon^2}$, coincide only if $d_\ell^2 = -\frac{h^2}{4} \sin^2(\frac{\ell h}{2}) = 0$, which is equivalent to $\frac{\ell h}{2} \in \pi \mathbb{Z}$. In most cases, this does not occur and therefore the case $g^2=\det(G)$ does not arise. Then, we obtain $g^2<\det(G)$ and the eigenvalues are distinct and satisfy $|\lambda_{1,2}| =1$, so that stability follows.
\end{remark}

\begin{remark}
\label{rem:x-dependence}
If the cross sections are space-dependent, the above analysis is not valid. In practice, the CFL condition is replaced by a worst-case condition. This means that we replace $\sigma_a$ and $\sigma_t$ in~\eqref{eq:condition_CFL} and \eqref{eq:condition_phi} by its maximum and minimum,
\begin{equation}
\sigma_{a,\max} = \max_{x} \sigma_a(x) \quad\text{and}\quad \sigma_{t,\min} = \min_{x} \sigma_t(x) \;,
\end{equation}
respectively.
\end{remark}

\begin{remark}
\label{rem:stability2D}
 In two dimensions, one can expect that the stability result from Proposition~\ref{prop:main_stability_result} carries over with the following changes. We replace $h= \min(\Delta x, \Delta y)$ and add a factor of $\frac12$ in front of the time step to account for the presence of growth rates in each of the two spatial dimensions.
\end{remark}

\section{Numerical Results}
\label{sec:numerics}
In this section, we consider different numerical test cases to demonstrate the performance of our scheme. Since we did not examine  boundary conditions, we only consider examples, where the solution is compactly supported away from the boundary. We implemented periodic boundary conditions, so that there is no influence of any discretization of boundary values.

The numerical calculations are performed using the two-dimensional scheme described in Section~\ref{sec:method} with the stability conditions from Section~\ref{sec:stability}. This means, we first choose the number of grid points ($N\times N$) for the staggered grids corresponding to the test case. Then, we determine the maximal time step (cf.\ Proposition~\ref{prop:main_stability_result}, Remark~\ref{rem:x-dependence}, and Remark~\ref{rem:stability2D})
\begin{equation}
\label{eq:numerics-CFL}
\Delta t := 0.9\; \cdot \tfrac12\;\min\left\{\tfrac{1}{\sigma_{a,\max}}, \max\{\tfrac{1}{2}\varepsilon h,\; \tfrac 1 4 h^2\sigma_{t,\min}\}\right\}
\end{equation}
and define the relaxation parameter
\begin{equation}
\label{eq:numerics-phi}
\phi := \begin{cases} h\frac{\sigma_{t,\min}}{2\varepsilon^3}\;, & h \sigma_t\leq 2\varepsilon\;,\\
\frac{1}{\varepsilon^2} &\text{otherwise}.\end{cases}
\end{equation}
with $h:= \frac1N$, $\sigma_{a,\max}:= \max_{x} \sigma_a(x)$, and $\sigma_{t,\min} := \min_{x} \sigma_t(x)$. The angular discretization uses a Gaussian quadrature with $16$ points on the interval $[0,1]$ for $\lambda$. As the quadrature points are mapped to the directions $\xi$ and $\eta$ with \eqref{eq:lambda}, we obtain $16$ points per quadrant. In all test cases, we compare the numerical solution on a grid where the parameter~$\varepsilon$ is resolved to a grid on which it is underresolved, thus demonstrating the AP property.

In the remainder of this section, we describe the test cases and the numerical results in detail. We consider four test cases to show different aspects of the AP property. First, we focus on the $\varepsilon$-dependence and investigate the convergence order in different regimes. In the second and third test cases, there are large spatial differences in the cross sections. The second test case is continuous and rotationally invariant, whereas in the third test case the material cross sections and  the source term are discontinuous. These two test cases intend to demonstrate the performance in multiscale problems. The last test case investigates the stability of the scheme-dependent on the choice of the relaxation parameter $\phi$.

\subsection{Convergence order}
\label{subsec:convergence_order}
We examine the order of convergence with respect to the spatial variable. We expect first or second order convergence depending on the used CFL condition. If a hyperbolic condition is used, the time step is proportional to $h$. As the explicit Euler method is used for the time discretization, we cannot expect more than first order convergence in $h$. Whereas if the parabolic condition is used, the time step is proportional to $h^2$. Then, the explicit Euler method predicts $\mathcal{O}(h^2)$ convergence. Moreover, centered differences, which are used for the spatial discretization, are as well a second order approximation in $h$. Thus, we expect that the error is proportional to $\mathcal{O}(h)$ when the hyperbolic condition is used, and $\mathcal{O}(h^2)$, respectively, when the parabolic condition is used. To estimate the convergence order, we compute the $\ell^2$-error $E(N)$ between the solution computed on a $N\times N$ grid and a reference solution. Using two different values $N_1$ and $N_2$, we then estimate the convergence order by
\begin{equation}
E_{N_1}^{N_2} = -\frac{\log(E(N_1))-\log(E(N_2))}{\log(N_1)-\log(N_2)}\;.
\end{equation}

\subsubsection{Method of manufactured solutions}
For the method of manufactured solutions (MMS), we first choose some function $f(t,x,y,\xi,\eta)$ and compute a corresponding source term and an initial condition, so that the chosen function is a solution of the transport equation. Let
\begin{equation}
f(t,x,y,\xi,\eta) = \exp(-t)\sin(2\pi x)^2\sin(2\pi y)^2(1+\eta^2)
\end{equation}
with $(x,y)\in [0,1]^2$. Further, let the scattering cross sections be given by $\sigma_a = 0$ and $\sigma_s = 1$. Then, the corresponding source term is given by
\begin{equation}
Q(t,x,y,\xi,\eta) = \partial_t f+\varepsilon v \cdot \nabla_x f - \frac{1}{\varepsilon^2} \left[ \frac{1}{2\pi} \int_\Omega f dv' \right],
\end{equation}
and the initial condition is given by
\begin{equation}
f(t=0,x,y,\xi,\eta) =  \sin(2\pi x)^2\sin(2\pi y)^2(1+\eta)^2.
\end{equation}
We use the source term and the initial condition to compute a solution with the above scheme. For different grid sizes and different values of $\varepsilon$, we compare the computed densities with the analytic density
\begin{equation}
\rho(t,x,y) = \frac{1}{2\pi} \int_\Omega f dv'
\end{equation}
at time $t=0.1$. The results are shown in Figure~\ref{fig:ex_mms} and Table~\ref{tab:ex_mms}. They confirm second order convergence in the parabolic case. In the hyperbolic case, the convergence order is even slightly higher than expected.

\noindent
\begin{minipage}{\textwidth} 
\begin{minipage}[b]{0.55\textwidth}
	\centering
	\includegraphics[width=\textwidth,trim = 0mm 0mm 15mm 0mm, clip=true]{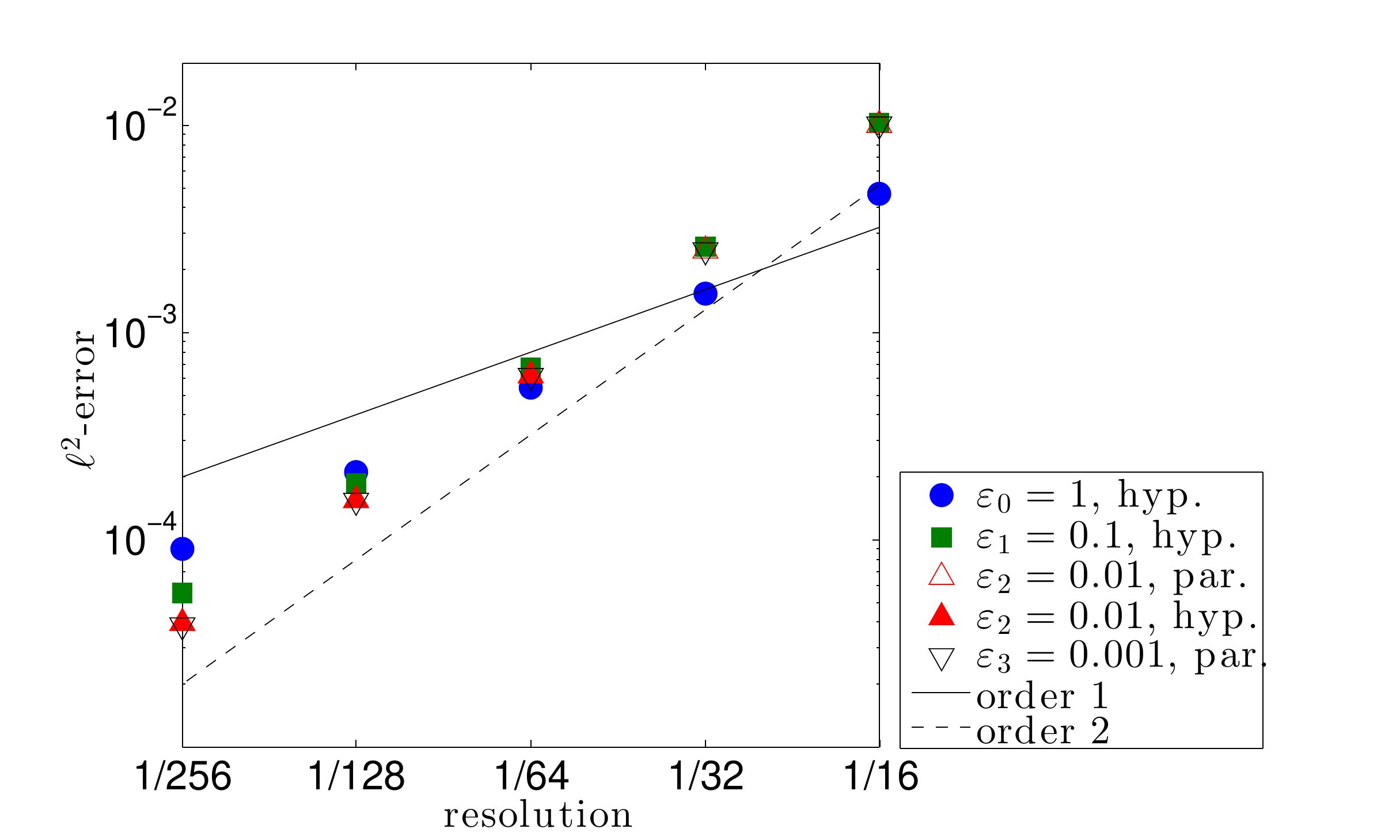}
	\captionof{figure}{Convergence order (MMS): $\ell^2$-error as a function of the spatial resolution. Hyperbolic (filled markers) or parabolic (empty markers) CFL condition.}
	\label{fig:ex_mms}
\end{minipage}\hfill
\begin{minipage}[b]{0.43\textwidth}
	\centering
	\begin{tabular}{ l | c c c c c c }
		\hline \footnotesize
 		&$E_{16}^{32}$& $E_{32}^{64}$ & $E_{64}^{128}$ & $E_{128}^{256}$\\
 		\hline
		$\varepsilon_0$      & 1.60 & 1.50 & 1.35 & 1.23 \\
		$\varepsilon_1$   & 1.98 & 1.93 & 1.86 & 1.76 \\
		\cdashline{2-3}
		$\varepsilon_2$ & 2.02  &\multicolumn{1}{:c:} { 2.00} & 1.99 & 1.97 \\
		\cdashline{3-5}
		$\varepsilon_3$& 2.02 & 2.01 & 2.00& 2.00 \\
		\hline
	\end{tabular}\\[2ex]
	\captionof{table}{Convergence order (MMS): The term $E_{N_1}^{N_2}$ is the convergence rate when going from $N_1\times N_1$ to $N_2\times N_2$ grid points for a fixed mean free path $\varepsilon_k= 10^{-k}$, $k=0,1,2,3$. The dashed line indicates the switch from the hyperbolic to the parabolic condition.}
	\label{tab:ex_mms}
\end{minipage}
\end{minipage}

\subsubsection{Gauss test}
We consider an example case with a smooth initial condition and isotropic scattering
\begin{equation}\label{eq:setup_ex_1}
\begin{aligned}
&f(t=0,x,y,v) = \tfrac{1}{4\pi \cdot 10^{-2}} \exp(-\tfrac{x^2+y^2}{4\cdot 10^{-2}})\quad \text{for} \quad   (x,y)\in [-1,1]\times[-1,1]\;,\\
& Q=0\;,\quad \sigma_t=\sigma_s=1\;,\quad \sigma_a=0\;,\quad\text{and}\quad \varepsilon = 1,\;10^{-1},\;10^{-2}\;. 
\end{aligned}
\end{equation}
Then, we compute the density $\rho$ at time $t=0.1$ for different grid sizes and different values of $\varepsilon$, so that the CFL condition~\eqref{eq:numerics-CFL} changes form hyperbolic to parabolic. As a reference solution, we use a highly resolved solution with $512\times512$ grid points. Table~\ref{tab:ex_1} and Figure~\ref{fig:ex_1} agree with the above assertion, showing first order convergence when the hyperbolic condition holds and second order, respectively, when the parabolic condition holds.

\noindent
\begin{minipage}{\textwidth} 
\begin{minipage}[b]{0.55\textwidth}
	\centering
	\includegraphics[width=\textwidth,trim = 0mm 0mm 15mm 0mm, clip=true]{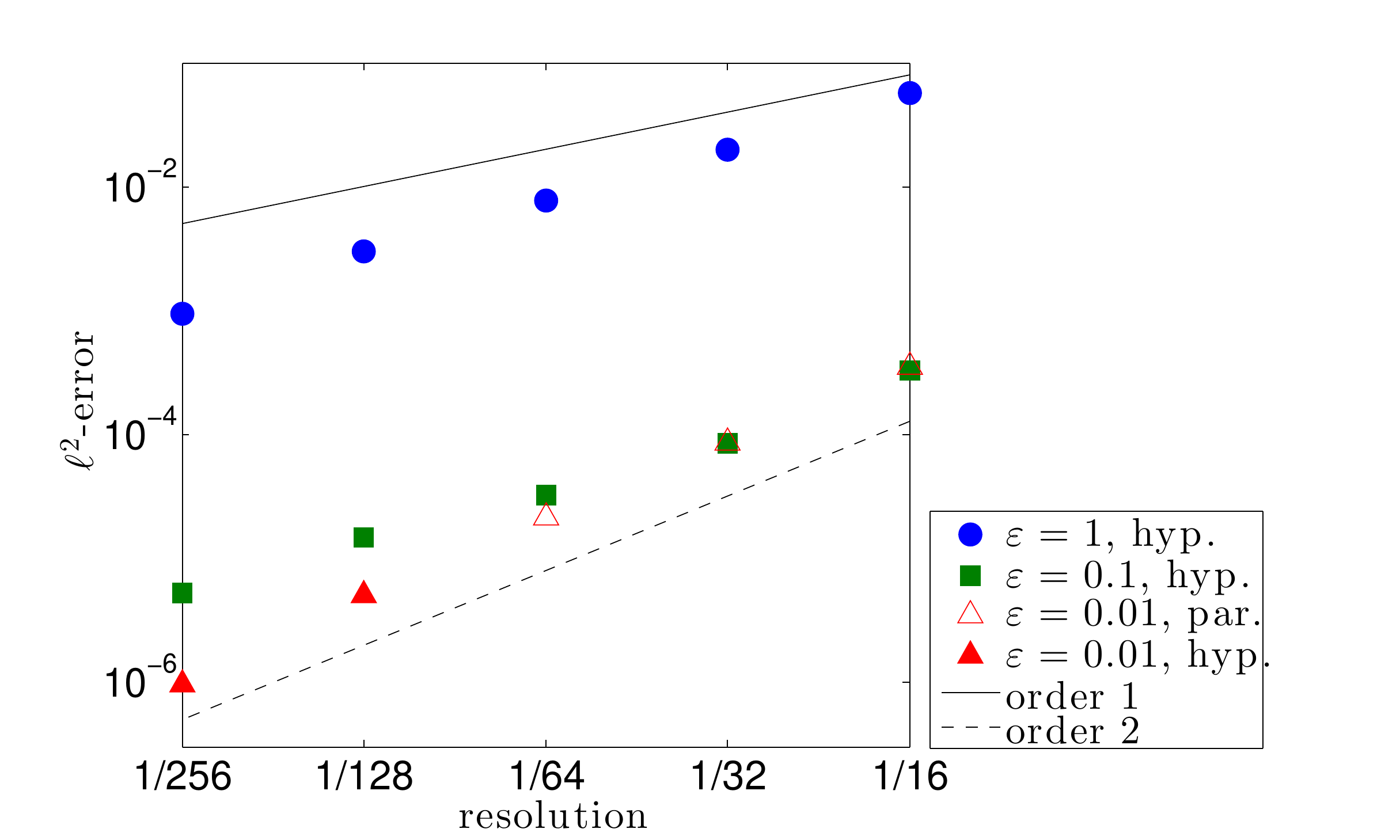}
	\captionof{figure}{Convergence order (Gauss test): $\ell^2$-error as a function of the spatial resolution. Hyperbolic (filled markers) or parabolic (empty markers) CFL condition.}
	\label{fig:ex_1}
\end{minipage}\hfill
\begin{minipage}[b]{0.43\textwidth}
	\centering
	\begin{tabular}{ l | c c c c }
		\hline\footnotesize
 		& $E_{16}^{32}$& $E_{32}^{64}$ & $E_{64}^{128}$ & $E_{128}^{256}$\\
	 	\hline
		$\varepsilon_0$       & 1.50 & 1.38 & 1.36 & 1.66 \\ 
		$\varepsilon_1$    & 1.96 & 1.37 & 1.15 & 1.50 \\ 
		\cdashline{2-4}
		$\varepsilon_2$  & 2.03 & 2.01 &\multicolumn{1}{:c:} { 2.09} & 2.41 \\ 
		\cdashline{4-5}
		\hline
	\end{tabular}\\[2ex]
	\captionof{table}{Convergence order (Gauss test): The term $E_{N_1}^{N_2}$ is the convergence rate when going from $N_1\times N_1$ to $N_2\times N_2$ grid points for a fixed mean free path $\varepsilon_k=10^{-k}$, $k=0,1,2$. The dashed line indicates the switch from the hyperbolic to the parabolic condition.}
\label{tab:ex_1}
\end{minipage}
\end{minipage}

\subsection{Variable scattering} 
In this test case, we examine the performance of the scheme, when the scattering is space-dependent. Compared to the previous test case, we fix the scaling parameter $\varepsilon$ and modify the scattering cross section. Let
\begin{equation}
\begin{aligned}
&f(t=0,x,y,v) = \tfrac{1}{4\pi \cdot 10^{-2}} \exp(-\tfrac{x^2+y^2}{4\cdot 10^{-2}})\quad \text{for} \quad (x,y) \in [-1,1]\times[-1,1]\;,\\
& \varepsilon=\tfrac{1}{100}\;,\quad Q=0\;,\quad \sigma_a=0\;, \quad \text{and} \\
&\sigma_t(x,y) =\sigma_s(x,y)=
\begin{cases} c^4(c+\sqrt{2})^2(c-\sqrt{2})^2\;, & c= \sqrt{x^2+y^2}<1\;,\\
1 &\text{otherwise.} \end{cases}
\end{aligned}
\end{equation}
Note that the total cross section $\sigma_t(x,y)$ can be periodically extended to a \mbox{$C^2$-function} and $\frac{\sigma_t(x,y)}{\varepsilon}$ ranges from $0$ to $100$. This wide range compared to the size of the domain causes strong variations of the solution, which are a challenge for numerical schemes.

\begin{figure}
\centering
\begin{subfigure}{0.32\textwidth}
	\includegraphics[trim=7mm 0mm 7mm 0mm, clip=true, width=\textwidth]{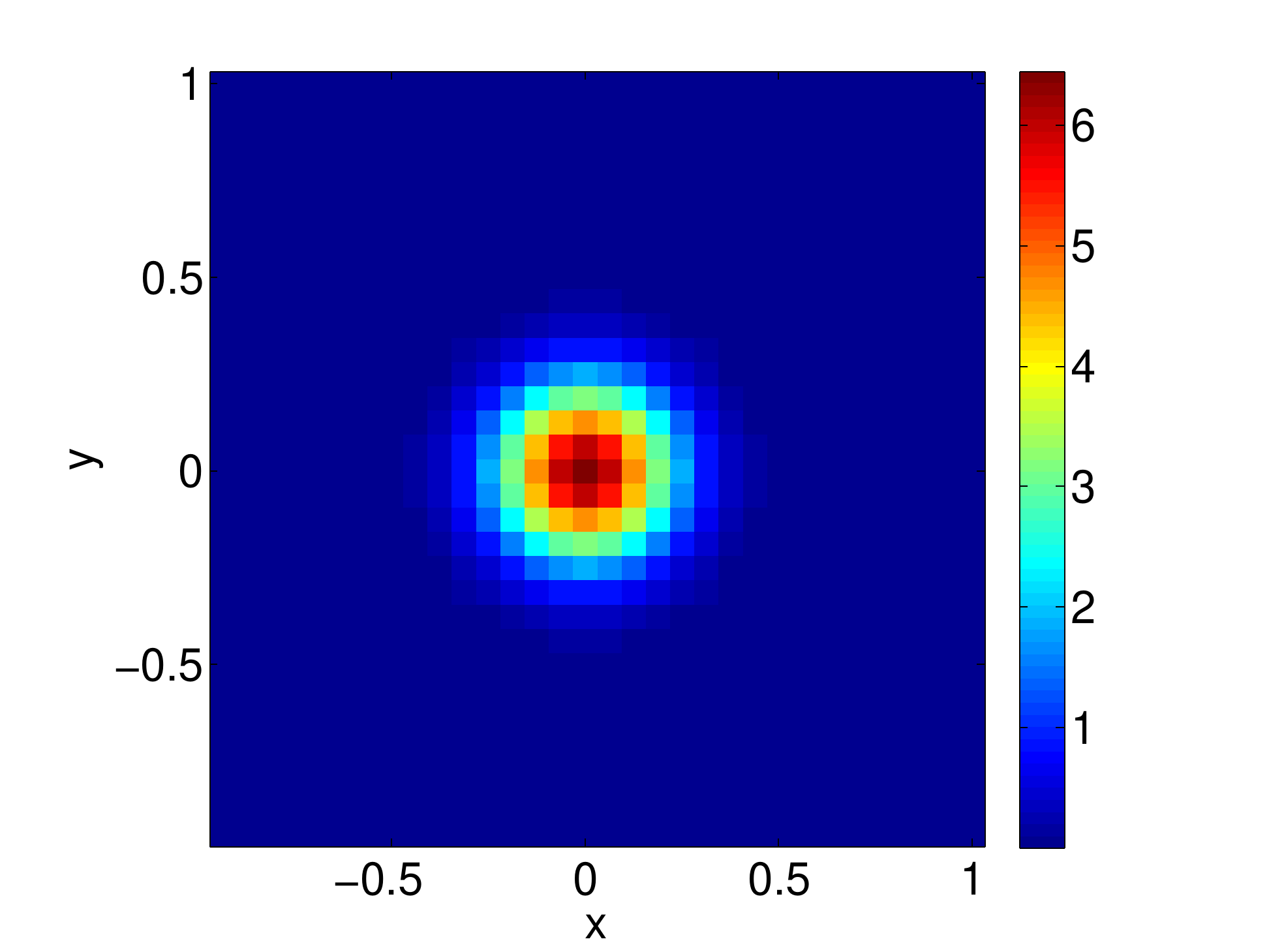}
	\includegraphics[trim=7mm 0mm 7mm 0mm, clip=true, width=\textwidth]{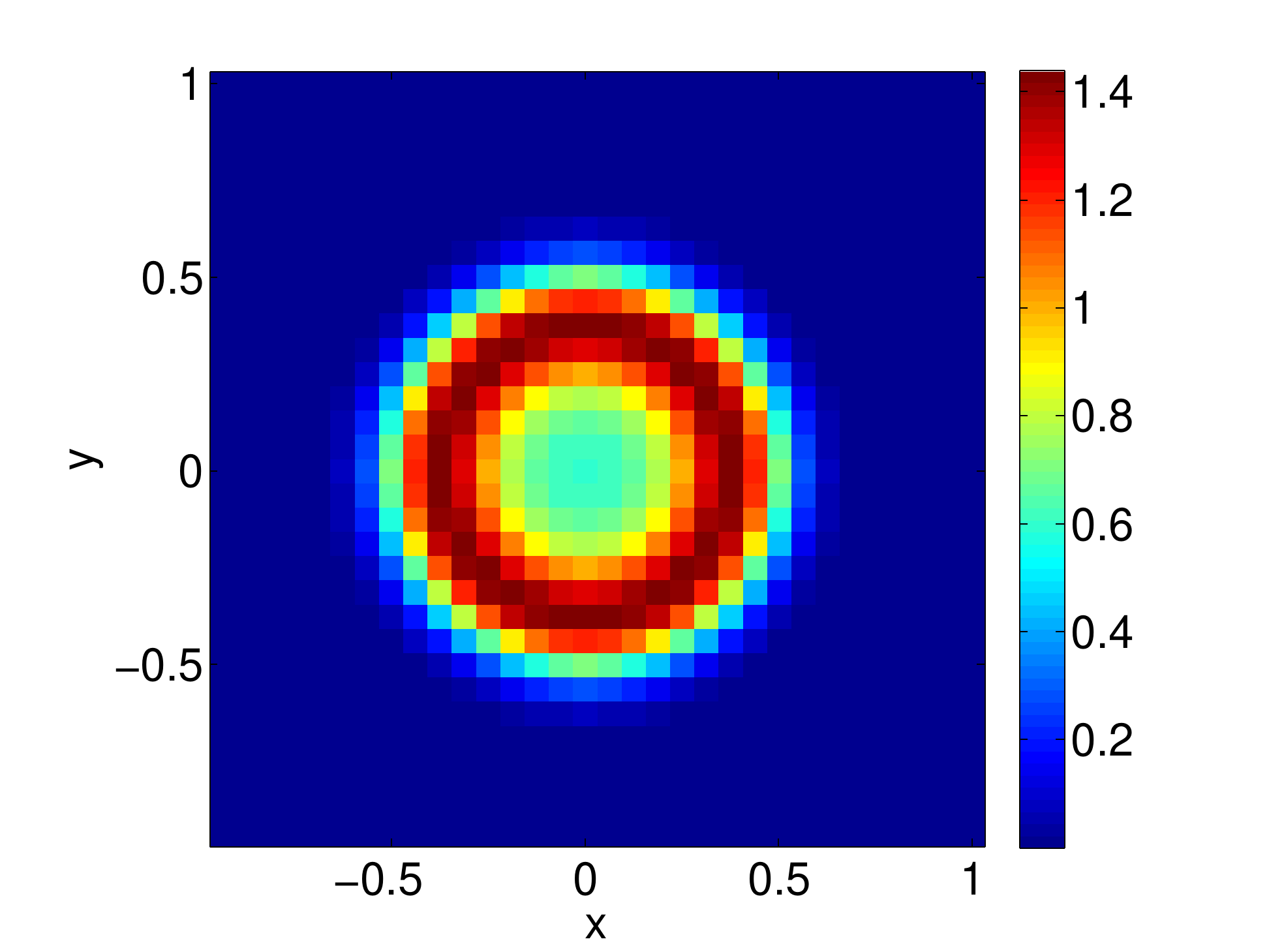}
	\includegraphics[trim=7mm 0mm 7mm 0mm, clip=true, width=\textwidth]{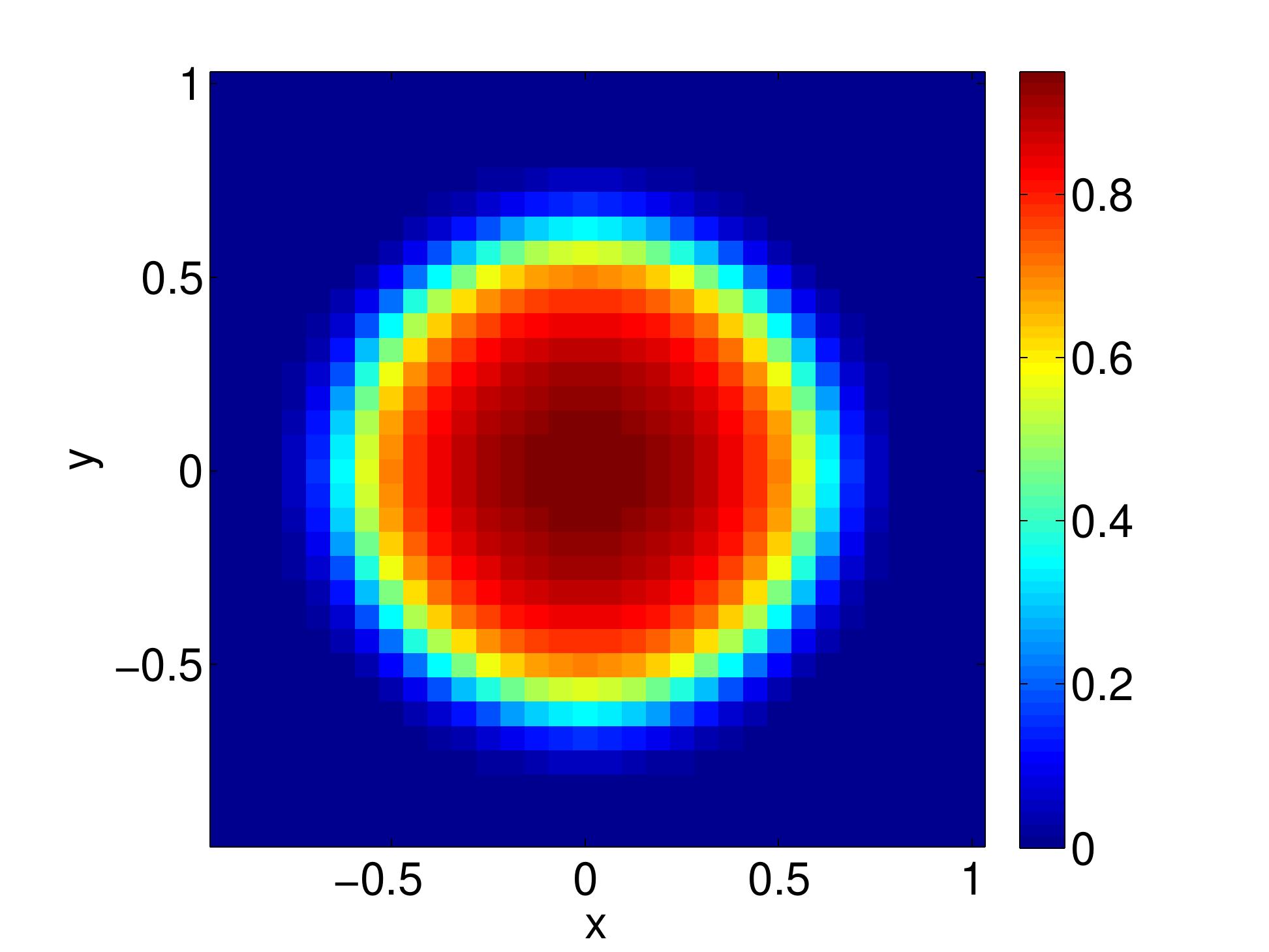}
\end{subfigure}
\begin{subfigure}{0.32\textwidth}
	\includegraphics[trim=7mm 0mm 7mm 0mm, clip=true, width=\textwidth]{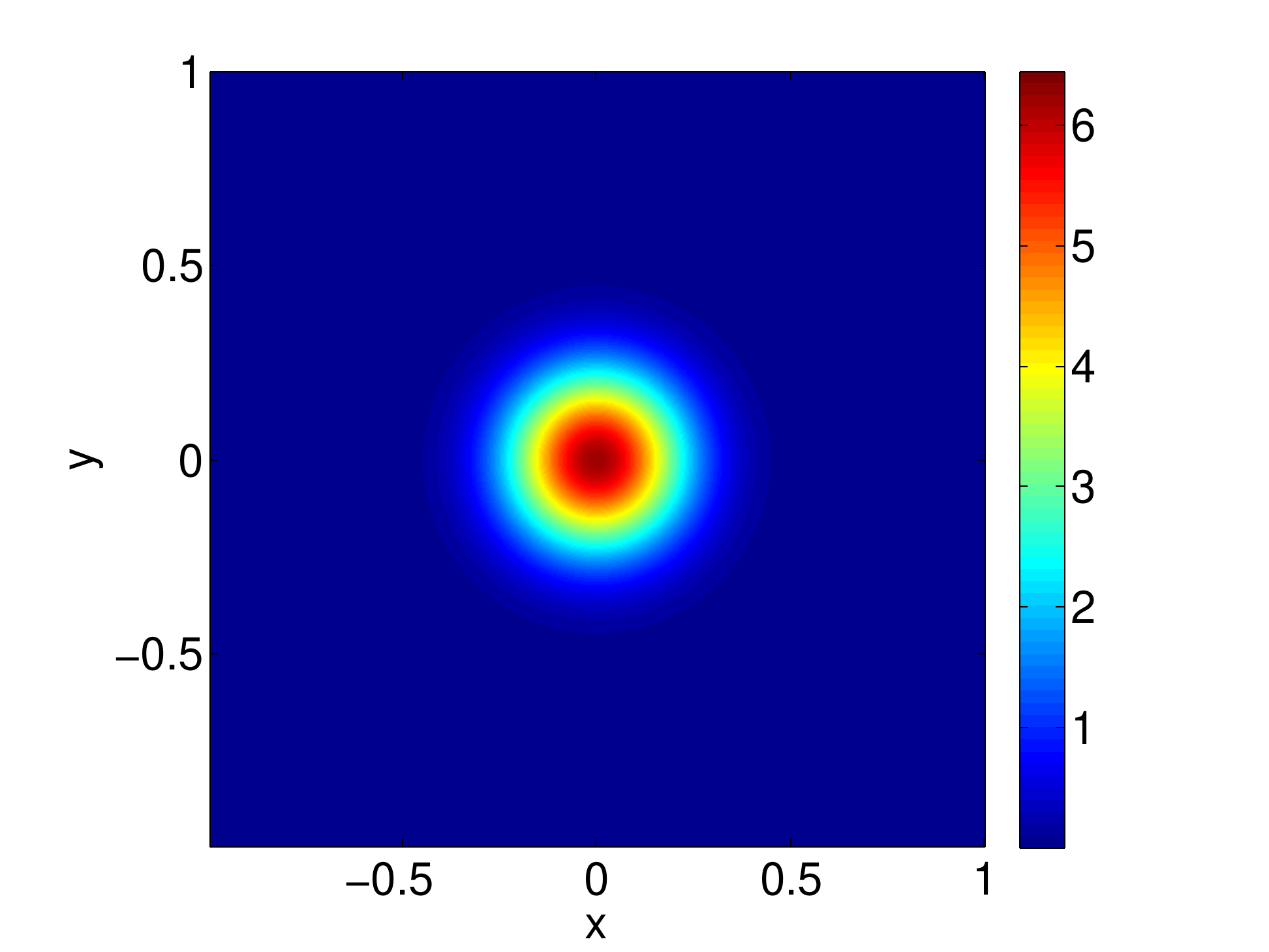}
	\includegraphics[trim=7mm 0mm 7mm 0mm, clip=true, width=\textwidth]{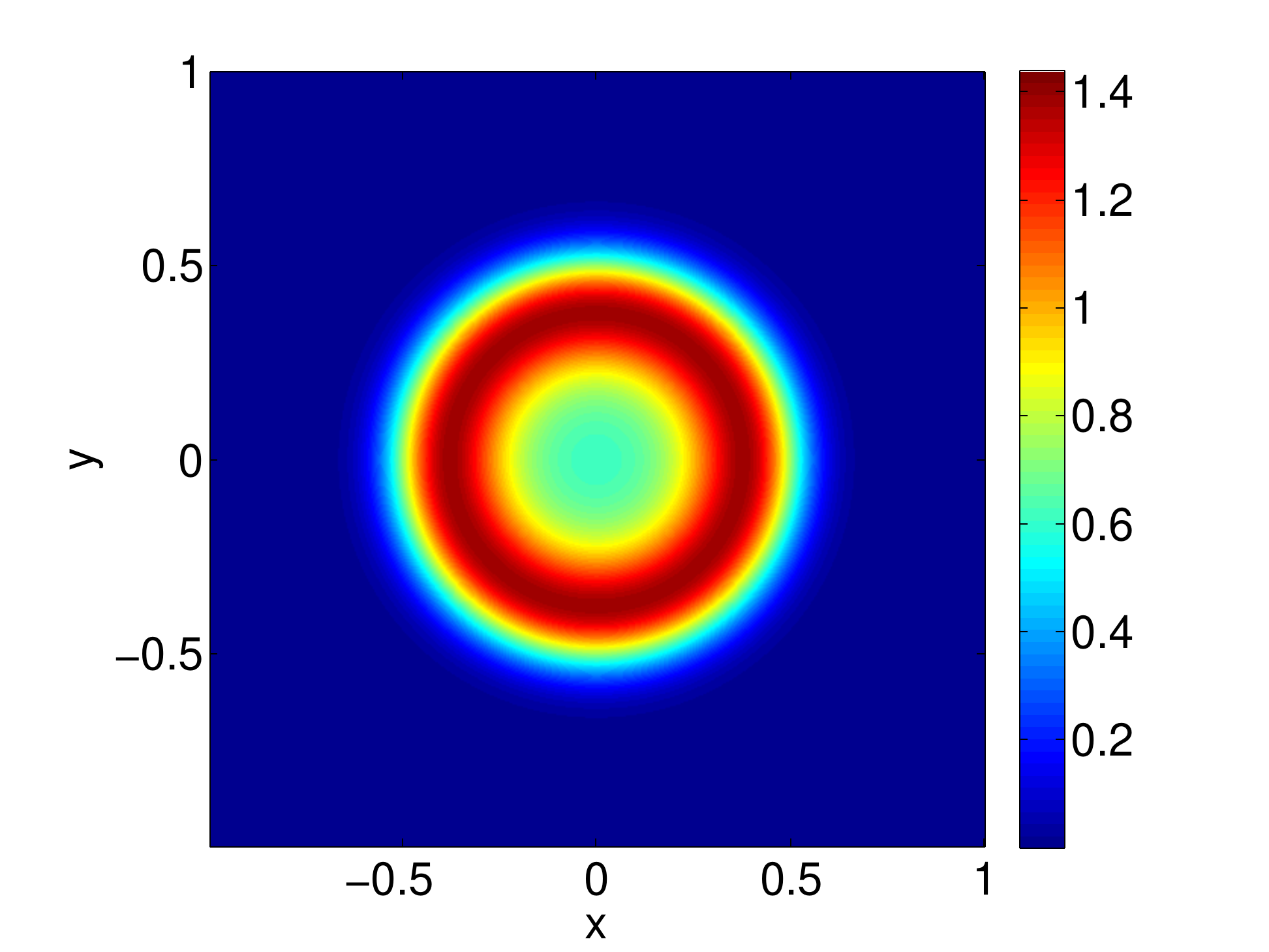}
	\includegraphics[trim=7mm 0mm 7mm 0mm, clip=true, width=\textwidth]{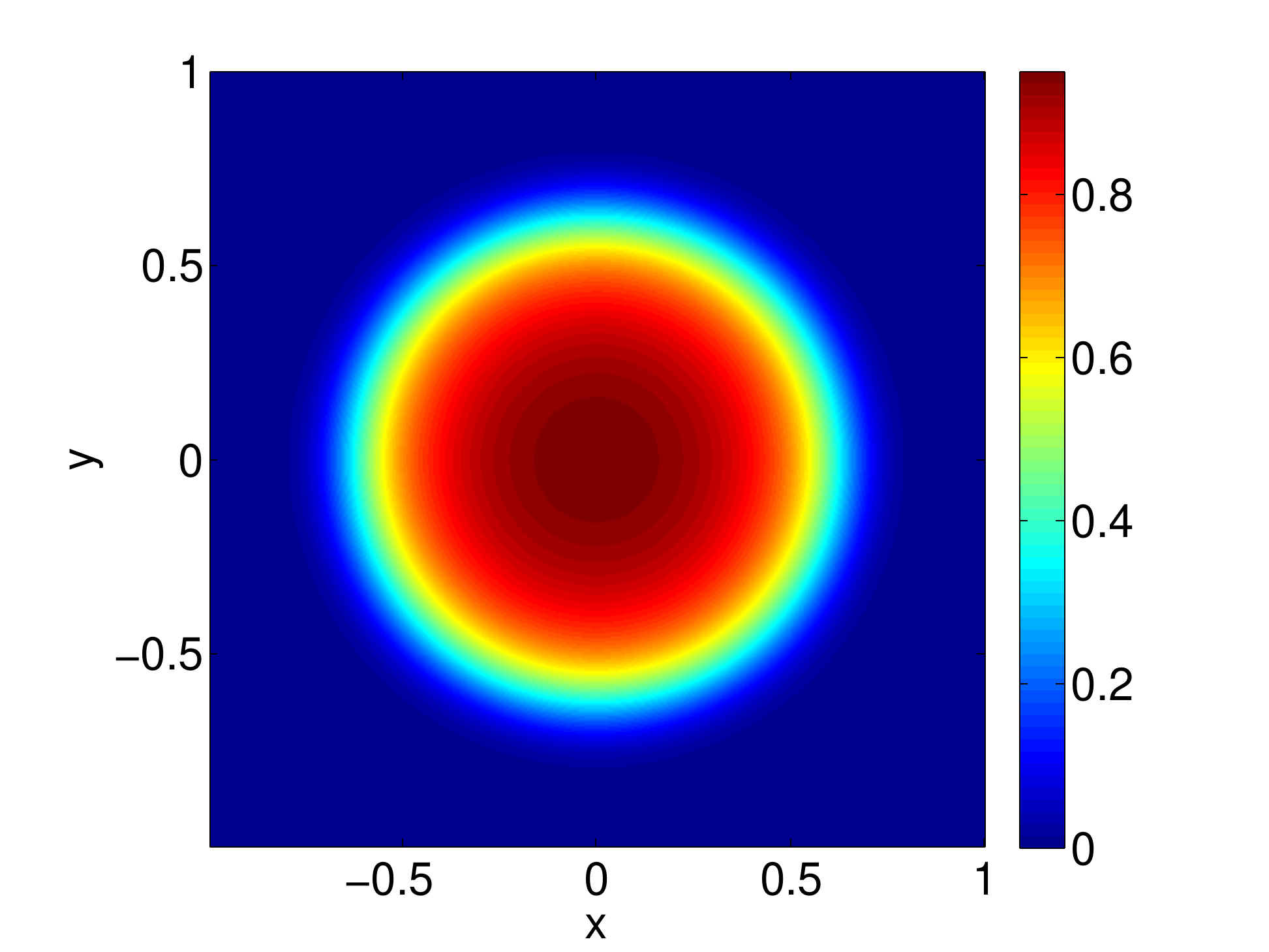}
\end{subfigure}
\begin{subfigure}{0.32\textwidth}
	\includegraphics[trim=7mm 0mm 7mm 0mm, clip=true, width=\textwidth]{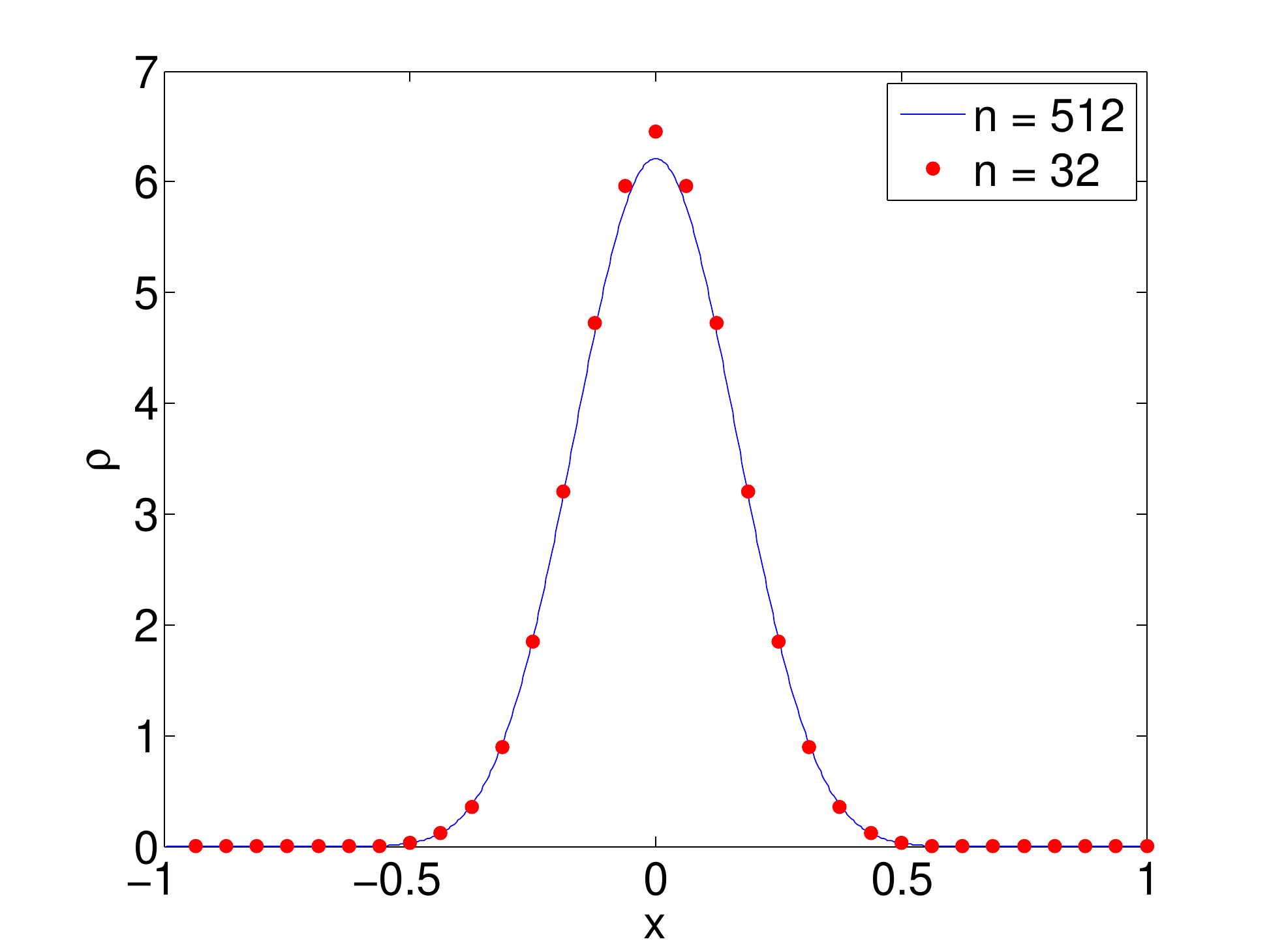}
	\includegraphics[trim=7mm 0mm 7mm 0mm, clip=true, width=\textwidth]{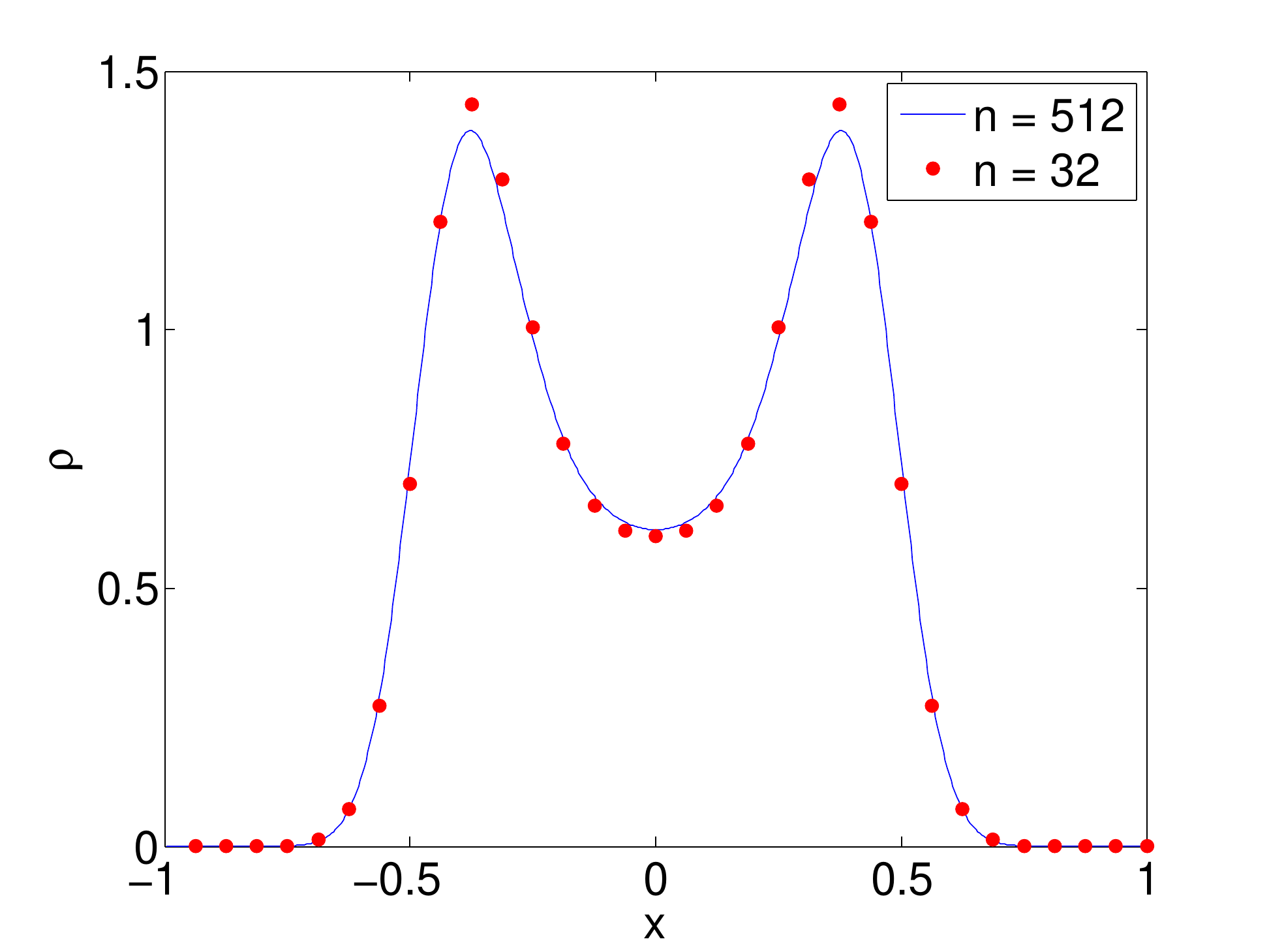}
	\includegraphics[trim=7mm 0mm 7mm 0mm, clip=true, width=\textwidth]{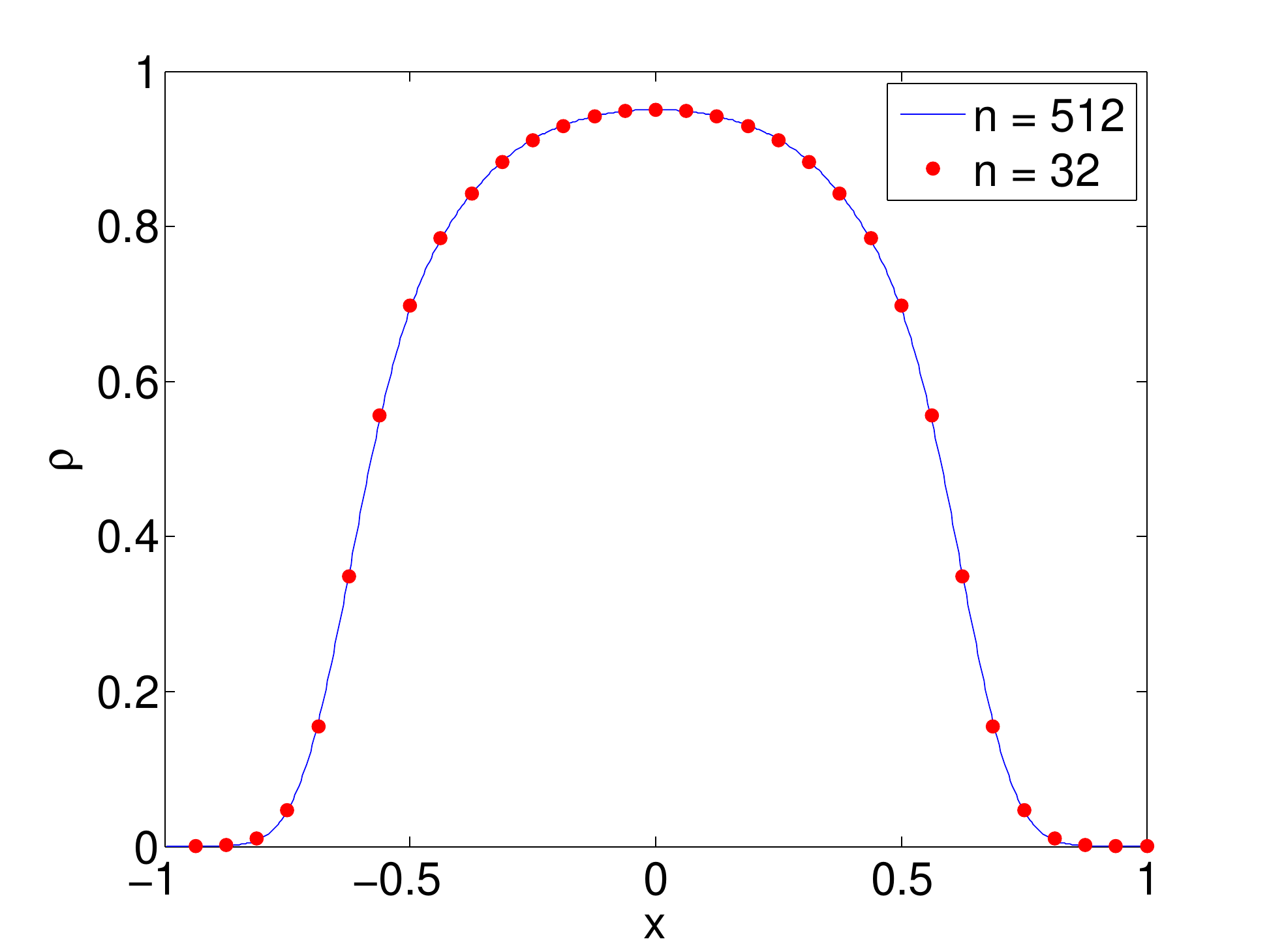}
\end{subfigure}
\caption{Variable scattering: Density at $\frac{t}{\varepsilon}=0.1$ (first row), $\frac{t}{\varepsilon}=0.5$ (second row), and $\frac{t}{\varepsilon}=1.0$ (third row), computed on a $32\times 32$ grid (first colum) or a $512\times512$ grid (second column). The third column shows the density on a cut along $y=0$.}
\label{fig:ex_2}
\end{figure}

We compute the solution up to time $t=\varepsilon$ on two different grids. One of the grids underresolves the length scale $\varepsilon=\tfrac{1}{100}$ ($32\times 32$ grid points) and the other one resolves it ($512\times 512$ grid points). Comparing the solution at different times ($t=\frac{1}{10}\varepsilon,\; \frac{1}{2}\varepsilon,\;\varepsilon$; see Figure~\ref{fig:ex_2}), we observe that the density, computed on the underresolved grid, matches the behavior of the density, computed on the resolved grid.\\

\subsection{Two material test} 
The two material test case is a slight modification of the lattice test, which was proposed in~\cite{Brunner-Holloway-2005}. It models a domain with different materials by discontinuous material cross sections and a discontinuous source term in space.

\begin{figure}
\begin{subfigure}{0.325\textwidth}
	\includegraphics[trim=15mm 0mm 15mm 0mm, clip=true, width=\textwidth]{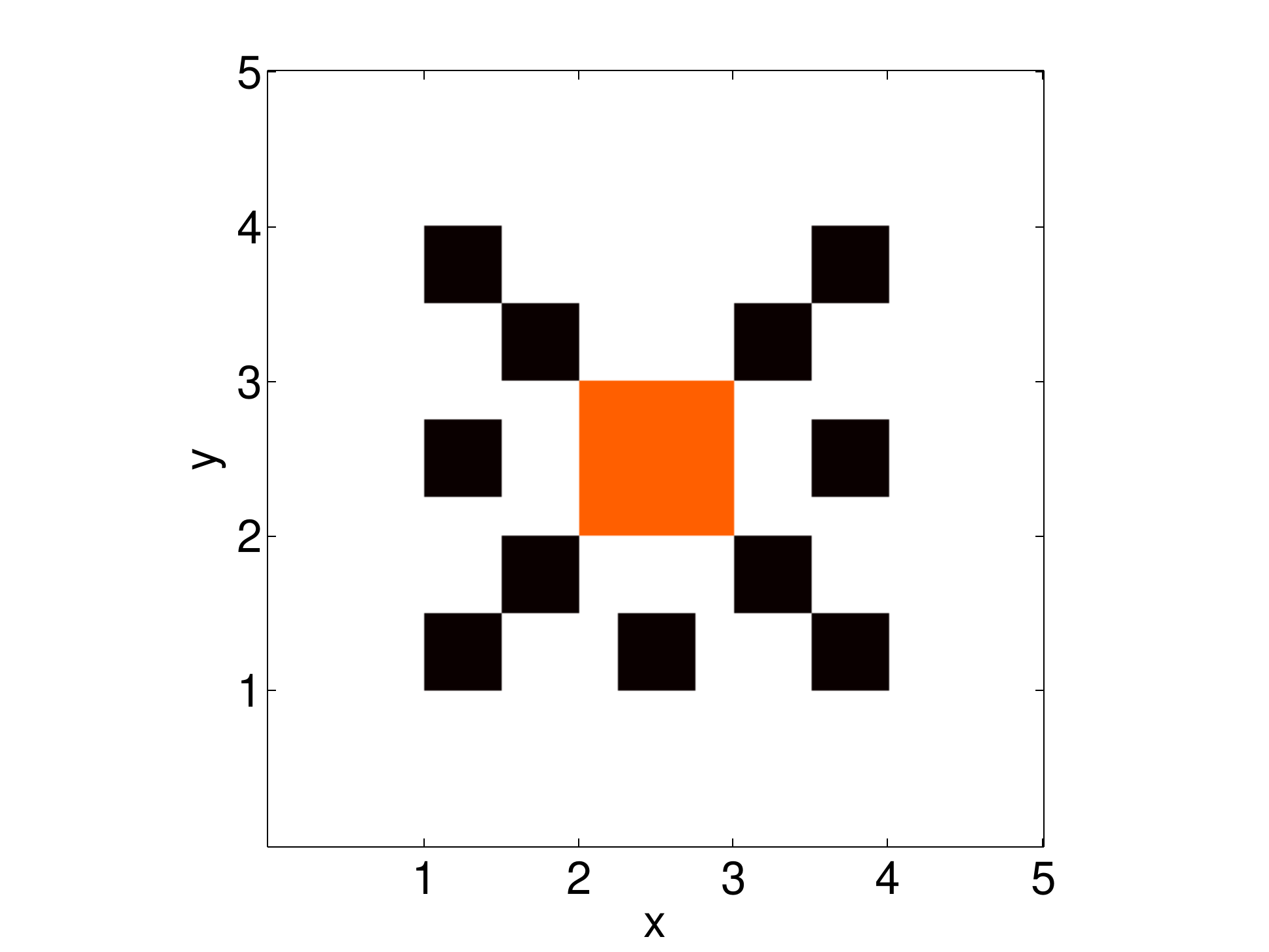}
	\caption{Geometry}
	\label{fig:ex_3_geometry}
\end{subfigure}\hfill
\begin{subfigure}{0.325\textwidth}
	\includegraphics[trim=15mm 0mm 15mm 0mm, clip=true, width=\textwidth]{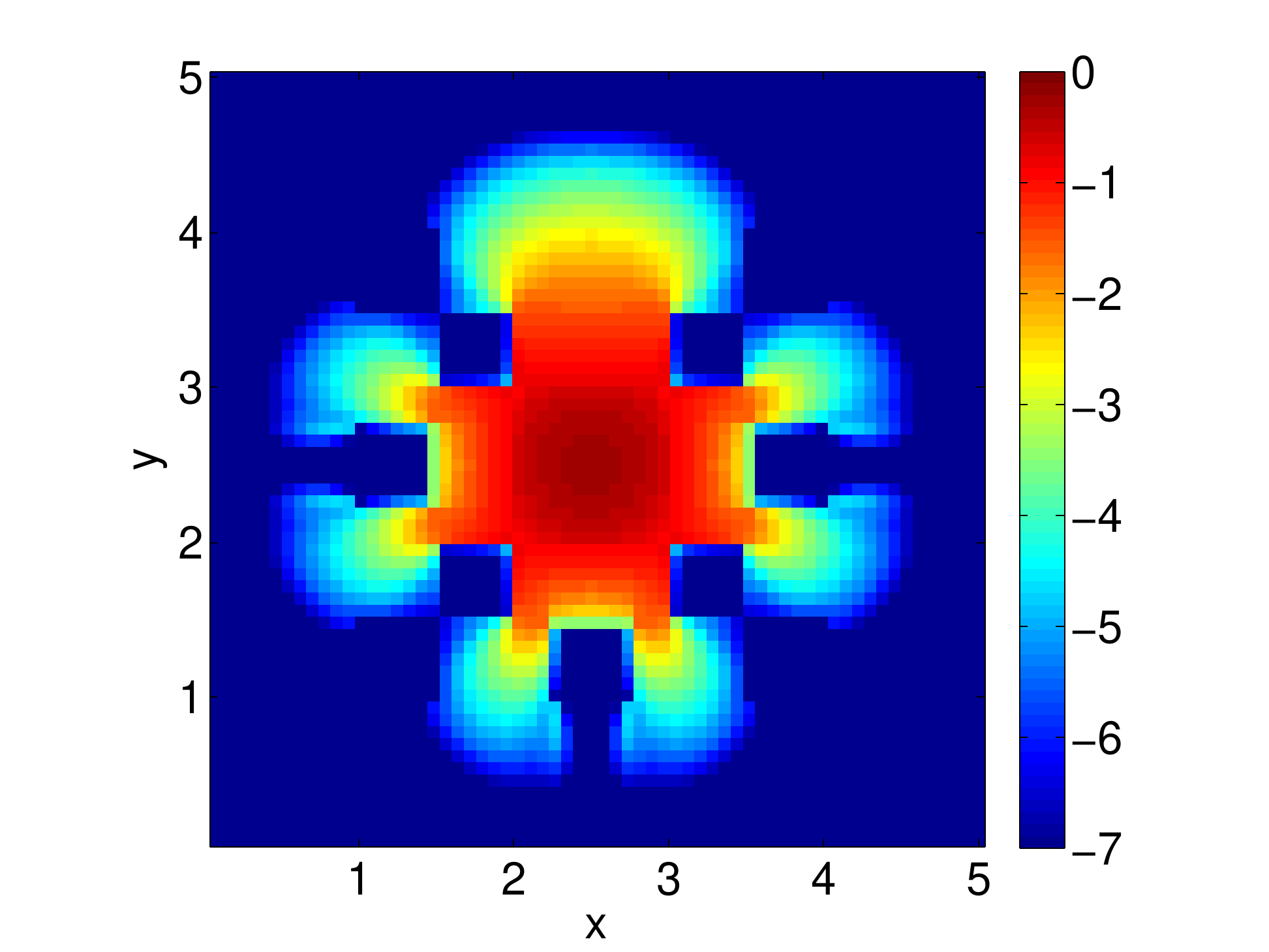}
	\caption{$64\times 64$ grid}
	\label{fig:ex_3_under-resolved}
\end{subfigure}\hfill
\begin{subfigure}{0.325\textwidth}
	\includegraphics[trim=15mm 0mm 15mm 0mm, clip=true, width=\textwidth]{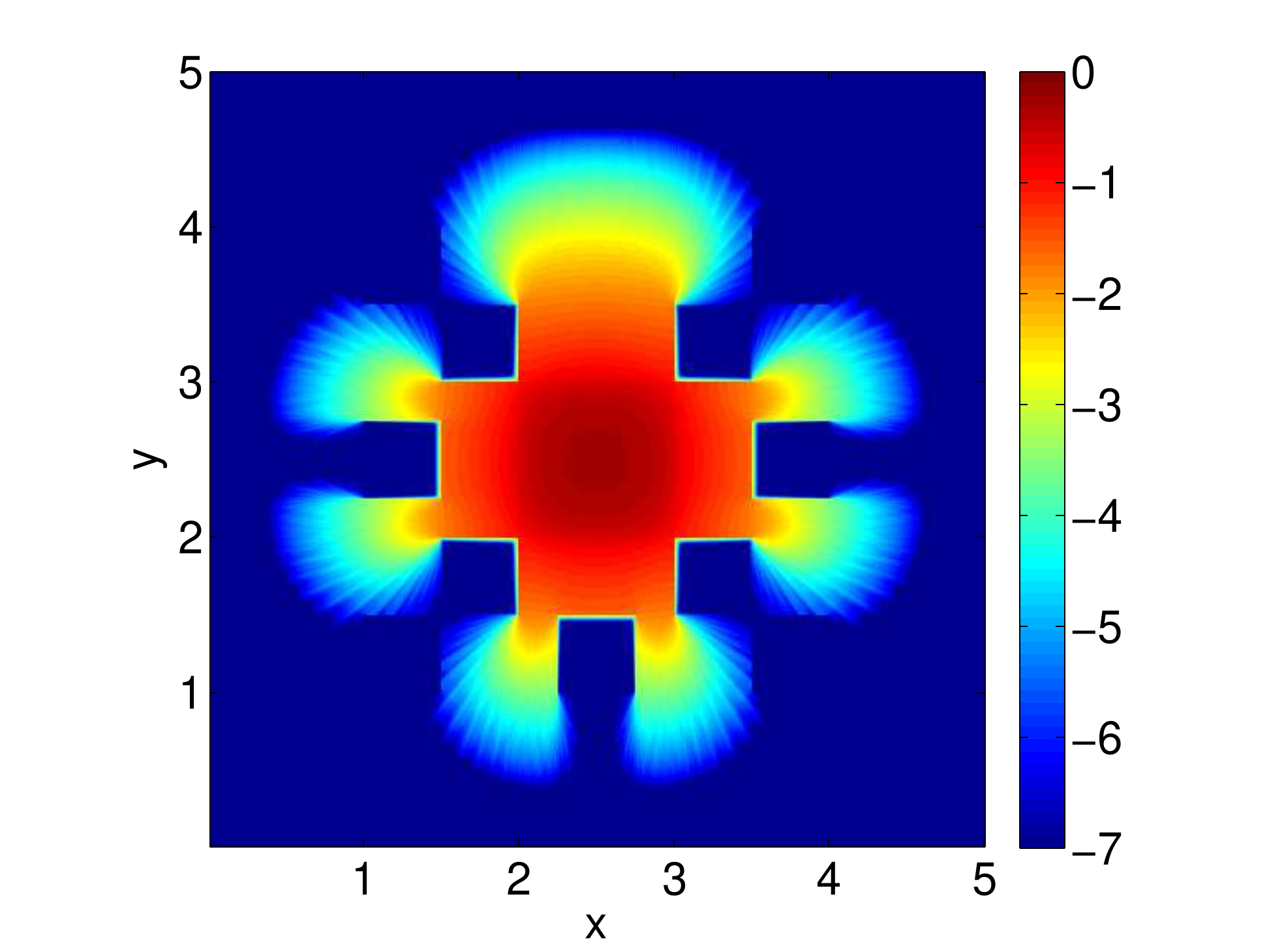}
	\caption{$512\times 512$ grid}
	\label{fig:ex_3_resolved}
\end{subfigure}
\caption{Two material test: (a) Geometry -- source (orange), purely scattering $\sigma_t=\sigma_s=1$ (white and orange), purely absorbing $\sigma_t=\sigma_a=100$ (black). (b) and (c) Density $\rho$ at $t=1.7$, computed on a $64\times 64$ grid (b) or $512\times512$ grid (c). Logarithmic scaling, values are limited to seven orders of magnitude.}
\label{fig:ex_3}
\end{figure}

In this problem, the computational domain is a $5\times 5$ square. Most of the domain is purely scattering, except for some purely absorbing squares of size $0.5$, which are distributed around an isotropic source in the middle of the domain
\begin{equation}
 \quad Q(x,y) = \begin{cases} 1 \;, &(x,y)\in[2,3]^2\;,\\ 0 & \text{otherwise.} \end{cases}
\end{equation}
In the absorbing spots (cf.\ Figure~\ref{fig:ex_3_geometry}), the absorption coefficient jumps from $0$ to $100$, while the scattering coefficient jumps from $1$ to $0$. Thus, there are diffusive and kinetic regimes, although the scaling parameter satisfies $\varepsilon=1$. We obtain a rapid change of the solution at the transition zones, which may cause difficulties in the numerics.

We compute the density up to time $t=1.7$ on a coarse grid ($64\times 64$) and on a fine grid ($512\times 512$). The solutions are shown in Figure~\ref{fig:ex_3}. Again, we observe that the solution on the underresolved grid resembles the solution on a grid that is resolved. In the case of the resolved solution, the oscillations near the beam edges are due to the angular discretization. They are the well-known ray effects for finite discrete velocity models (cf.~\cite{Brunner-2002} and references therein, as well as~\cite{Lathrop-1968,Miller-Reed-1977}).

\subsection{Relaxation parameters and stability} 
In our final test, we consider different relaxation parameters. Proposition~\ref{prop:main_stability_result} suggests an upper bound on the relaxation parameter $\phi$, which in the hyperbolic case is more restrictive than in the parabolic case. We expect that in certain example cases our scheme becomes unstable if $\phi$ is too large.

\begin{figure}
\begin{subfigure}{0.49\textwidth}
	\includegraphics[width = \textwidth]{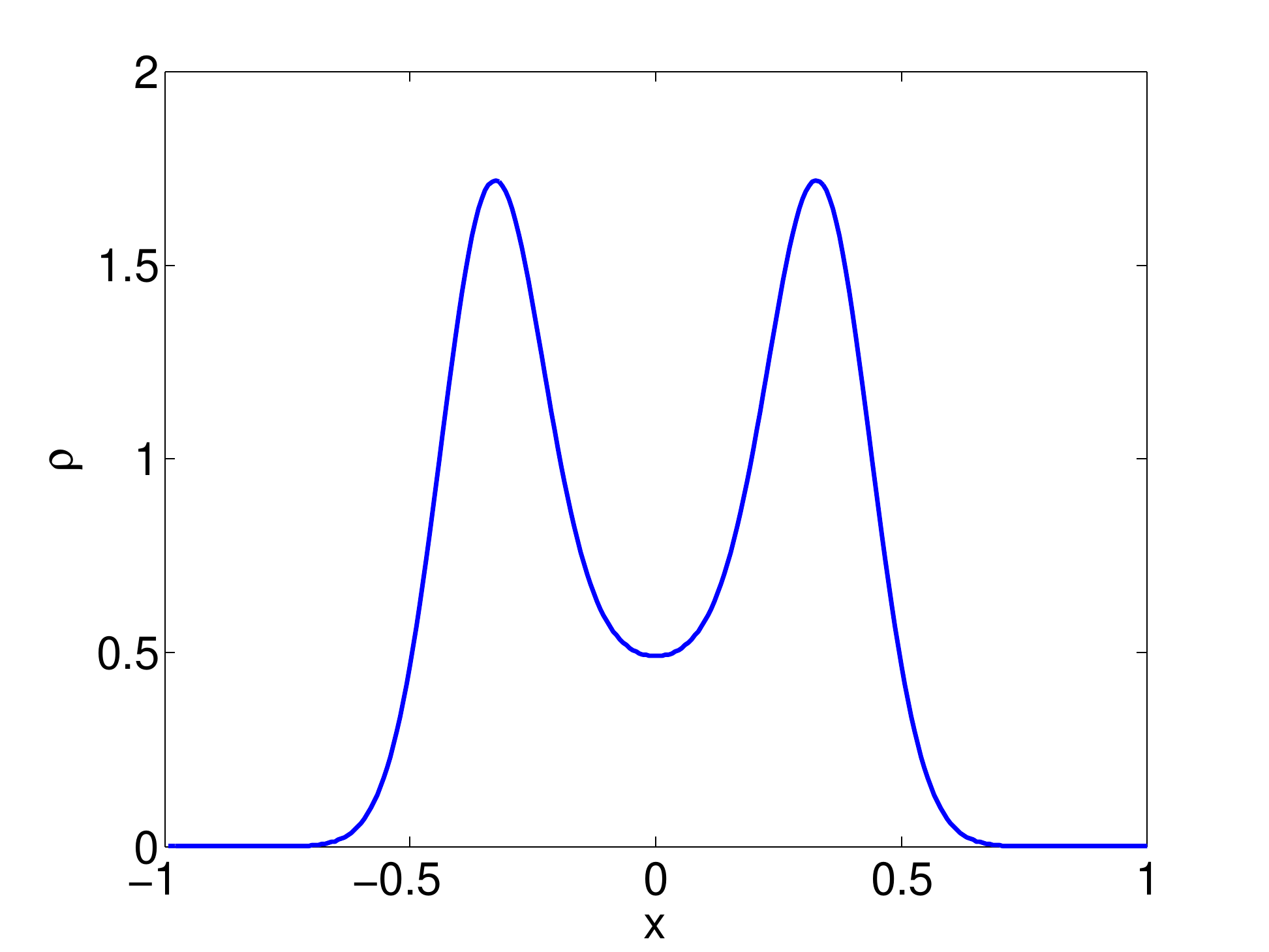}
	\caption{$\phi = \tfrac{h\sigma_t}{2\varepsilon^3}$}
	\label{fig:ex_4_stable}
\end{subfigure}
\begin{subfigure}{0.49\textwidth}
	\includegraphics[width = \textwidth]{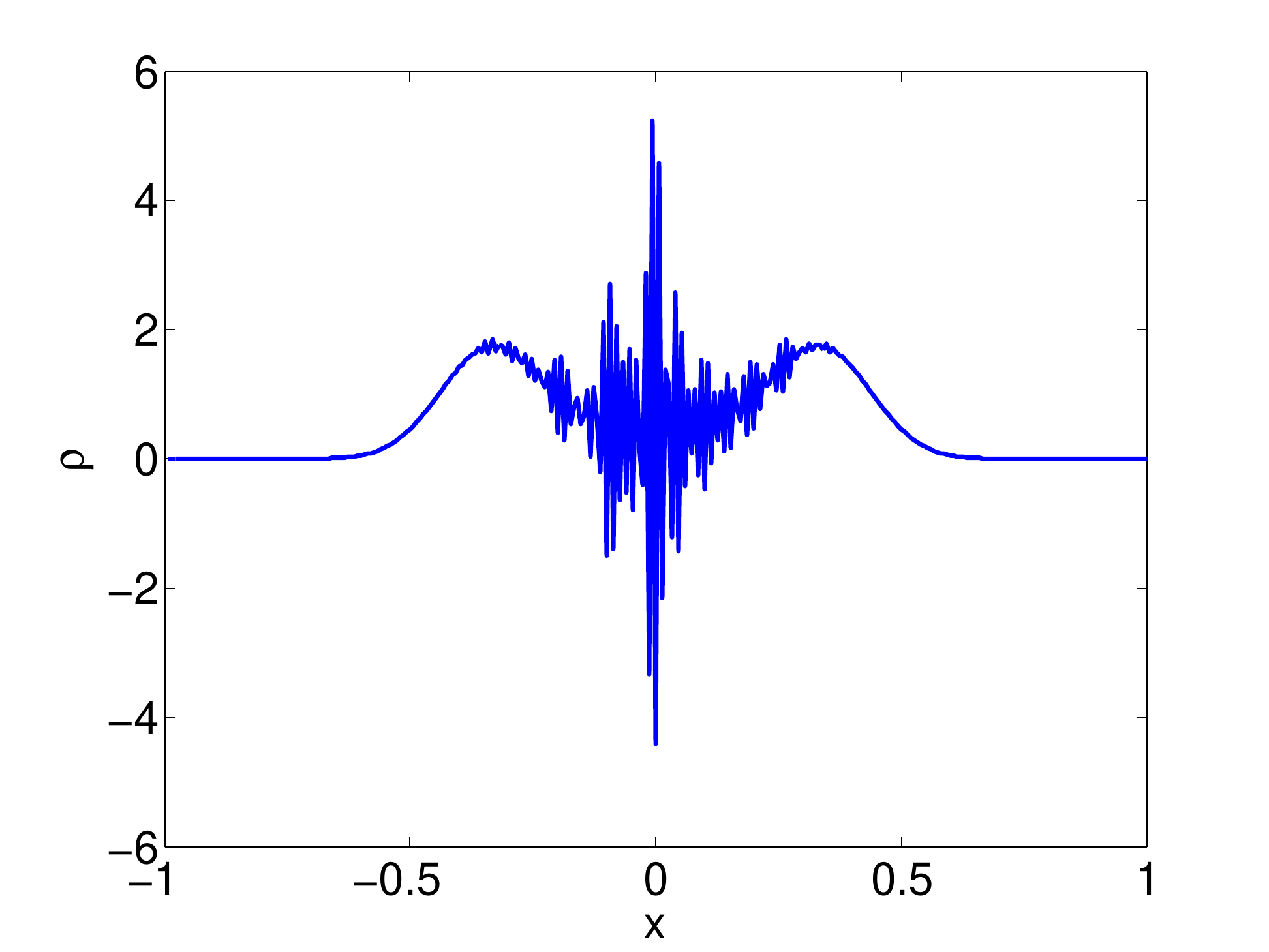}
	\caption{$\phi = \tfrac{1}{\varepsilon^2}$}
	\label{fig:ex_4_blow-up}
\end{subfigure}
\caption{Stability: Density on a cut along $y=0$, computed on a $300 \times 300$ grid up to time $t=0.36$ using different values of $\phi$.}
\label{fig:ex_4}
\end{figure}

Similar to the Gauss test, let
\begin{equation}
\begin{aligned}
&f(t=0,x,y,v) = \tfrac{1}{4\pi \cdot 5\times 10^{-3}} \exp(-\tfrac{x^2+y^2}{4\cdot 5\times 10^{-3}}) \quad\text{for}\quad (x,y)\in [-1,1]\times[-1,1]\;,\\
& \varepsilon=1\;,\quad Q=0\;,\quad \sigma_t=\sigma_s=1\;,\quad \sigma_a=0\;,\quad  N=300\;,\quad\text{and}\quad t=0.36\;.
\end{aligned}
\end{equation}
Then, we compute the density on a $N\times N= 300\times 300$ grid up to time $t=0.36$ using different relaxation parameters
\begin{equation}
\phi_1 = \tfrac{h\sigma_t}{2\varepsilon^3} = \tfrac{10}{3}\times 10^{-4} \quad\text{and} \quad \phi_2 =\tfrac1{\varepsilon^2} = 1.
\end{equation}
In the first case, the relaxation parameter satisfies the assumption of Proposition~\ref{prop:main_stability_result} and the solution is stable (see Figure~\ref{fig:ex_4_stable}), whereas in the second case, the assumption is violated and the solution starts to blow up (see Figure~\ref{fig:ex_4_blow-up}). As a consequence, the upper bound on the relaxation parameter in Proposition~\ref{prop:main_stability_result} can in general not be substituted by the less restrictive upper bound $\phi\leq \tfrac{1}{\varepsilon^2}$.

\section{Extensions}
\label{sec:ext}
\subsection{Boundary conditions}
In the numerical examples, we consider only periodic boundary conditions. However, one can use different approaches to implement, for instance, inflow boundary conditions. In the following, we discuss how to adapt two different methods~\cite{Klar-1998,Jin-Pareschi-Toscani-2000} to our scheme. Let the computational domain be given by $\mathbf{X}=[0,L_x]\times[0,L_y]$ for some $L_x,\;L_y>0$. Then, inflow boundary conditions on a continuous level are given by
\begin{align}
f(t,\mathbf{x},v) = f^{\text{in}}(t,\mathbf{x},v) \quad \text{for} \quad n\cdot\mathbf{x}>0,\; x\in\partial \mathbf{X},
\end{align}
where $n$ is the outer normal vector and $f^{\text{in}}$ describes the inflow.

\begin{figure}[htb]
\centering
\begin{tikzpicture} [x=1.5cm,y=1.25cm]
    \foreach \i in {0,1,2} {
        \draw [black, dashed] (\i,-1/2) -- (\i,2.5);
        \draw [black, dashed] (-1/2,\i) -- (2.5,\i) ;
    }
    \foreach \i in {0,1,2,3} {
        \draw [black, dashed] (\i-1/2,-1/2) -- (\i-1/2,2.5) ;
        \draw [black, dashed] (-1/2,\i-1/2) -- (2.5,\i-1/2);
    }
     \draw [black, thick] (0,0) -- (0,2) ;
     \draw [black, thick] (0,0) -- (2,0) ;
     \draw [black, thick] (2,2) -- (0,2) ;
     \draw [black, thick] (2,2) -- (2,0) ;
\node[draw,circle,inner sep=2.5pt,fill,red] at (1,1) {};

  \foreach \x in {1,2}{
      \foreach \y in {1,2}{
        \node[draw,circle,inner sep=2.5pt,fill,red] at (\x-1/2,\y-1/2) {};
      }
    }
    
\foreach \y in {1,2}{
	\node[draw,diamond,inner sep = 2.pt,fill,blue] at (1,\y-1/2) {};
	}
\foreach \x in {1,...,2}{
	\node[draw,diamond,inner sep = 2.pt,fill,blue] at (\x-1/2,1) {};
	}

\foreach \y in{0,2}{
	\node[draw,circle,inner sep=2.5pt,red] at (1,\y) {};
	\node[draw,circle,inner sep=2.5pt,red] at (\y,1) {};
	}
\foreach \y in{1,2}{
 \foreach \x in{0,2}{
	\node[draw,diamond,inner sep=2.5pt,blue] at (\x,\y-1/2) {};
	\node[draw,diamond,inner sep=2.5pt,blue] at (\y-1/2,\x) {};
	}
}
\foreach \x in{0,3}{
\foreach \y in{1,2}{
	\node[draw,circle,inner sep=2.5pt,red] at (\x-1/2,\y-1/2) {};
	\node[draw,circle,inner sep=2.5pt,red] at (\y-1/2,\x-1/2) {};
	}
}
\foreach \y in{0,3}{
	\node[draw,diamond,inner sep=2.5pt,blue] at (1,\y-1/2) {};
	\node[draw,diamond,inner sep=2.5pt,blue] at (\y-1/2,1) {};
}
    \node [below] at (0,-.5) {$0$};
    \node [below] at (2,-.5) {$L_x$};
    \node [left] at (-.5,0) {$0$};
    \node [left] at (-.5,2) {$L_y$};
\end{tikzpicture}
\caption{Black solid line denotes the boundary of the domain $\partial \mathbf{X}$; filled markers denote inner points; empty markers denote boundary and ghost points.}\label{fig:ghost-points}
\end{figure}
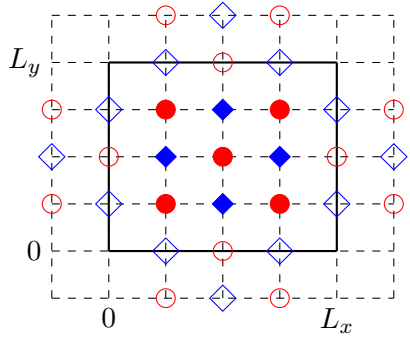

Similar to \cite{Klar-1998,Golse-Klar-1995,Klar-1998-II}, the solution of a kinetic half-space problem, describing the zeroth order diffusion boundary conditions, can be approximated to determine the outflow function $f^{\text{out}}(t,\mathbf{x},v)$:
\begin{align}
f(t,\mathbf{x},v) = f^{\text{out}}(t,\mathbf{x},v) \quad \text{for} \quad n\cdot\mathbf{x}<0,\; x\in\partial \mathbf{X}.
\end{align}
This determines completely the density $f$ on the boundary and with this also the parities on the boundary. Considering, for instance, the left boundary $x =0$, the boundary conditions of the parities $\pr{1}$, $\pj{1}$ are given by
\begin{equation} \label{eq:boundary-r1-j1}
\begin{aligned}
\pr{1}(\xi,\eta)+\varepsilon\pj{1}(\xi,\eta) &= f(\xi,-\eta) =f^{\text{in}}(\xi,-\eta)\;,\\
\pr{1}(\xi,\eta)-\varepsilon\pj{1}(\xi,\eta) &= f(-\xi,\eta) =f^{\text{out}}(-\xi,\eta)\;.
\end{aligned}
\end{equation}
In one dimension, the parity $\pr{1}$ can be placed on the grid $i\Delta x$, $i=0,\ldots,N_x$ with $N_x=\frac{L_x}{\Delta x}$ and $\pj{1}$ on the staggered grid $(i+\frac{1}{2})\Delta x$, $i=-1,\ldots,N_x$. Then $\pr{1}$ lies on the boundary and $\pj{1}$ can be interpolated using the ghost points $-\frac12$ and $N_x+\frac12$ \cite{Klar-1998}. This can be extended to two dimensions if we use a one-dimensional interpolation orthogonal to the boundary and add suitable ghost points (see Figure \ref{fig:ghost-points}). Then, the boundary condition~\eqref{eq:boundary-r1-j1} becomes for $\pr{1}$ and $j = 1,\ldots,N_y-1$ and $N_y = \frac{L_y}{\Delta_y}$
\begin{equation}
\pr{1}_{0,j}+\frac{\varepsilon}{2}\left(\pj{1}_{-\tfrac12,j} +\pj{1}_{\tfrac12,j} \right)=f^{\text{in}}_{0,j}\quad \text{and}\quad 
\pr{1}_{0,j}-\frac{\varepsilon}{2}\left(\pj{1}_{-\tfrac12,j} +\pj{1}_{\tfrac12,j} \right)=f^{\text{out}}_{0,j} 
\end{equation}
and for $\pj{1}$ and $j = \frac12,1+\frac12,\ldots,N-\frac12$
\begin{equation}
\frac{1}{2}\left(\pr{1}_{-\frac12,j} + \pr{1}_{\frac12,j}\right) +\varepsilon \pj{1}_{0,j} =f^{\text{in}}_{0,j} \quad \text{and}\quad 
\frac{1}{2}\left(\pr{1}_{-\frac12,j} + \pr{1}_{\frac12,j}\right) -\varepsilon \pj{1}_{0,j} =f^{\text{out}}_{0,j} \;.
\end{equation}
Note that the edges of the boundary can be dropped, because these points are needed neither for computing the inner points nor for the interpolation (due to the one-dimensional interpolation rule). This choice yields a system with the same number of equations and unknowns.

Alternatively, the boundary conditions presented in~\cite{Jin-Pareschi-Toscani-2000} can be used. The key idea is to approximate the j-unknowns on the boundary by
\begin{equation}
\pj{1} = -\xi \partial_x \pr{1} - \eta \partial_y \pr{1} \quad \text{and} \quad 
\pj{2} = -\xi \partial_x \pr{2} - \eta \partial_y \pr{2}
\end{equation}
and insert this into the inflow boundary parts. For example, at the left boundary $x=0$, the boundary condition for $\pr{1}$ is given by
\begin{equation}
f^{\text{in}}(\xi,-\eta) = \pr{1}(\xi,\eta)+\varepsilon\pj{1}(\xi,\eta) = \pr{1}(\xi,\eta)-\varepsilon \left(\xi \partial_x \pr{1} + \eta \partial_y \pr{1}\right) \;.
\end{equation}
To implement the boundary conditions in our scheme the spatial derivatives can, for instance, be approximated by one-sided finite-differences. There are multiple choices to do this, which need to be tested.

\subsection{Anisotropic scattering}
The simplicity of the scheme relies on the simplicity of the scattering operator, which is isotropic. In the following, we briefly consider the case of anisotropic scattering. Let the scattering operator be given by
\begin{align}
K f(v) = \int s(v,v') f(v') dv'
\end{align}
for some symmetric scattering kernel $s$, which depends on the cosine of the scattering angle $v\cdot v'$.

To compute the parity equations, the transport equation is rewritten for $v=(\xi,\eta),$ $(\xi,-\eta),$ $(-\xi,\eta)$, and $(-\xi,-\eta)$ with nonnegative $\xi,\eta\geq 0$ and added in correspondence with the definition of the parities~\eqref{eq:definition-parities}. This leads to a summation of the scattering operator of the form $Kf(v) \pm Kf(-v)$, which has to be expressed in terms of the parities. Since the scattering kernel satisfies $s(-v,v') = s(v,-v')$, we obtain
\begin{multline}\label{eq:SumScatteringOperator}
Kf(v) \pm Kf(-v) = \int s(v,v')f(v') dv' \pm \int  s(-v,v') f(v') dv' \\
= \int s(v,v')f(v') dv' \pm \int  s(v,v') f(-v') dv'
= \int s(v,v')(f(v')\pm f(-v') )dv' \;. 
\end{multline}
More precisely, we consider for $\xi,\;\eta\geq 0$
\begin{align}
\frac12 (Kf(\xi,-\eta)+ Kf(-\xi,\eta)) &= \int s(\xi,\eta,\xi',\eta') \frac12 \left( f(\xi',-\eta')+f(-\xi',\eta')\right) d(\xi',\eta')\;.
\end{align}
In this term, we recover the parities $\pr{1}$ and $\pr{2}$
\begin{align}
\frac12 ( f(\xi',-\eta')+f(-\xi',\eta')) = \begin{cases} \pr{1} & \text{for}\; \xi',\;\eta'\geq0\; \text{or}\; \xi',\; \eta'\leq 0\;,\\
\pr{2} & \text{otherwise.} \end{cases}\;.
\end{align}
Splitting the integral into the four quadrants and changing to positive $\xi',\eta'$, we obtain
\begin{multline}
\frac12 Kf(\xi,-\eta)+ Kf(-\xi,\eta)  =
\int_{\xi,\eta>0} (s(\xi,\eta,\xi',\eta')+ s(\xi,\eta,-\xi',-\eta')) \pr{1}(\xi',\eta') d(\xi',\eta') \\
+\int_{\xi,\eta>0}  (s(\xi,\eta,\xi',-\eta')+ s(\xi,\eta,-\xi',\eta')) \pr{2}(\xi',\eta') d(\xi',\eta')\;.
\end{multline}
Similarly, the other summations can be expressed in the parities. To simplify the notation, we define the following operators:
\begin{equation}
\begin{aligned}
K^+f(\xi,\eta) &:=\int_{\xi,\eta>0} (s(\xi,\eta,\xi',\eta')+ s(\xi,\eta,-\xi',-\eta')) f(\xi',\eta') d(\xi',\eta')\;, \\
K^-f(\xi,\eta) &:= \int_{\xi,\eta>0}  (s(\xi,\eta,\xi',-\eta')+ s(\xi,\eta,-\xi',\eta')) f(\xi',\eta') d(\xi',\eta')\;.
\end{aligned}
\end{equation}
Then, the parity equations for anisotropic  scattering are given by
\begin{equation}
\begin{aligned}
\partial_t \pr{1} + \xi \partial_x \pj{1} -\eta \partial_y \pj{1} 
&= -\frac{\sigma_s}{\varepsilon^2} ( \pr{1} - K^+\pr{1}-K^-\pr{2})-\sigma_a \pr{1} +Q\;,\\
\partial_t \pr{2} + \xi \partial_x \pj{2} +\eta \partial_y \pj{2} 
&= -\frac{\sigma_s}{\varepsilon^2} (\pr{2} -  K^+\pr{2}-K^-\pr{1})-\sigma_a \pr{2} +Q\;,\\
\partial_t \pj{1} + \frac\xi{\varepsilon^2} \partial_x \pr{1} -\frac\eta{\varepsilon^2} \partial_y \pr{1} 
&= -\frac{\sigma_s}{\varepsilon^2-} (\pj{1}-  K^+\pj{1}-K^-\pj{2})-\sigma_a \pj{1}\;,\\
\partial_t \pj{2} + \frac\xi{\varepsilon^2} \partial_x \pr{2} +\frac\eta{\varepsilon^2} \partial_y \pr{2} 
&= -\frac{\sigma_s}{\varepsilon^2} (\pj{2}-  K^+\pj{2}-K^-\pj{1}) -\sigma_a  \pj{2}\;.
\end{aligned}
\end{equation}
One could now apply a similar splitting as above, but in contrast to the isotropic case the implementation of the relaxation step is not straightforward. It is necessary to invert the integral operators at a potentially expensive computational cost (cf.~\cite{Jin-Pareschi-2000} for more details).

\section{Conclusions}
In this paper, we have introduced a two-dimensional AP scheme for the linear transport equation. The linear transport equation has the diffusion equation as an analytic asymptotic limit. For AP schemes the discretization has to be chosen such that the analytic limit is preserved at a discrete level and the scheme is uniformly stable with respect to the mean free path. Here, we used a parity-based time discretization combined with a staggered-grid spatial discretization.

We have shown that the spatial discretization has the desired AP
property. In particular, due to the use of staggered grids, a compact five-point stencil can be achieved in the limiting discrete diffusion limit. Furthermore, the parity-based time discretization is suitable for the use of staggered grids, as the coupling between the even and odd parities reduces the number of the required unknowns. In addition, we have presented a rigorous stability analysis for the same scheme in one dimension. This provides a condition on the relaxation parameter and a CFL condition. Finally, we have performed several numerical tests for the two-dimensional scheme, which demonstrate the AP property. Since staggered grids can easily be extended to three dimensions, there is a straightforward generalization of our method to three spatial dimensions. Although we did not test the method, we expect that it has similar properties.

In the future, it would be worthwhile to investigate the time discretization. Since our method uses a simple time-integration method (explicit Euler method), the convergence order is in general limited to one. To maintain a second order scheme, one could use some higher order IMEX 
time-integration method.  Another possible scope of future work is to apply
staggered grids in combination with a parity-based time discretization to other kinetic equations.

\bibliographystyle{abbrv}
\bibliography{ap}

\end{document}